\documentclass[english]{smfart}

\usepackage[francais,english]{babel}
\usepackage{amssymb, amscd, smfthm}
\usepackage{bull}

\newtheorem{mainth}{Theorem}

\newcommand{\Z}{\mathbb{Z}}
\newcommand{\Q}{\mathbb{Q}}
\newcommand{\F}{\mathbb{F}}
\newcommand{\D}{\mathcal{D}}
\newcommand{\Order}{\mathcal{O}}
\newcommand{\maxid}{\mathfrak{p}}
\newcommand{\onto}{\twoheadrightarrow}
\newcommand{\into}{\hookrightarrow}
\newcommand{\isomto}{\overset{\sim}{\to}}
\newcommand{\tensor}{\otimes}
\newcommand{\ctensor}{\mathbin{\hat{\otimes}}}
\newcommand{\dsum}{\oplus}
\newcommand{\Dsum}{\bigoplus}
\newcommand{\algcl}[1]{\overline{#1}}
\newcommand{\Ga}{\mathbf{G}_{a}}
\newcommand{\fGa}{\widehat{\mathbf{G}}_{a}}
\newcommand{\Gm}{\mathbf{G}_{m}}
\newcommand{\ab}{\mathrm{ab}}
\newcommand{\ur}{\mathrm{ur}}
\newcommand{\sep}{\mathrm{sep}}
\newcommand{\ram}{\mathrm{ram}}
\newcommand{\et}{\mathrm{et}}
\newcommand{\fpqc}{\mathrm{fpqc}}
\newcommand{\pr}{\mathrm{pr}}
\newcommand{\id}{\operatorname{id}}
\newcommand{\Gal}{\operatorname{Gal}}
\newcommand{\Hom}{\operatorname{Hom}}
\newcommand{\End}{\operatorname{End}}
\newcommand{\Ext}{\operatorname{Ext}}
\newcommand{\Ker}{\operatorname{Ker}}
\newcommand{\dlog}{\operatorname{dlog}}
\newcommand{\Spec}{\operatorname{Spec}}
\newcommand{\Isom}{\operatorname{Isom}}
\newcommand{\alg}[1]{\mathbf{#1}}
\newcommand{\AJ}[1]{\mathrm{\Xi}_{#1}}
\newcommand{\Proj}{\mathbb{P}}
\newcommand{\Supp}[1]{\mathrm{Supp}(#1)}
\newcommand{\sgpalg}[2]{#1 \langle #2 \rangle}
\newcommand{\pfpqc}[1]{(\mathrm{Perf} / #1)_{\mathrm{fpqc}}}
\newcommand{\set}[2]{
	\left\{
		#1 \;\; \vrule \;\; #2
	\right\}
}
\author{Takashi Suzuki}
\address{
	Department of Mathematics, University of Chicago,
	5734 S University Ave,
	Chicago, IL 60637, USA
}
\email{suzuki@math.uchicago.edu}
\title[Some remarks on the LCFT of Serre and Hazewinkel]
 {Some remarks \\ on the local class field theory \\ of Serre and Hazewinkel}
\alttitle{Quelques remarques sur la th\'eorie du corps de classes local de Serre et Hazewinkel}

\begin{document}

\frontmatter

\begin{abstract}
	We give a new approach for the local class field theory of Serre and Hazewinkel.
	We also discuss two-dimensional local class field theory in this framework.
\end{abstract}
\begin{altabstract}
	Nous donnons une nouvelle approche
	pour la th\'eorie du corps de classes local de Serre et Hazewinkel.
	Nous discutons \'egalement
	la th\'eorie du corps de classes local de dimension deux dans ce cadre.
\end{altabstract}

\altkeywords{Th\'eorie du corps de classes local, K-th\'eorie}

\maketitle


\mainmatter

\section{Introduction}

The purpose of this paper is twofold.
First, we give a geometric description of
the local class field theory of Serre and Hazewinkel (\cite{Ser61}, \cite[Appendice]{DG70})
in the equal characteristic case.
The main result is Theorem \ref{Thm:Main} below.
Second, we discuss its two-dimensional analog
with the aim to seek an analog of Lubin-Tate theory for two-dimensional local fields.
The discussion is taken place in Section \ref{Sect:TwoDim} with a partial result.

We formulate Theorem \ref{Thm:Main}.
The precise definitions of the terms used below are given in Section \ref{Sect:Prelim}.
Let $k$ be a perfect field of characteristic $p > 0$ and set $K = k((T))$.
The group of units of $K$ can be viewed as a proalgebraic group
(more precisely, a pro-quasi-algebraic group) over $k$
in the sense of Serre (\cite{Ser60});
we denote this proalgebraic group by $\alg{U}_{K}$.
For each perfect $k$-algebra $R$
(perfect means that the $p$-th power map is an isomorphism),
we have the group of its $R$-rational points
	\[
			\alg{U}_{K}(R)
		=
			\set{
				\sum_{i = 0}^{\infty}
					a_{i} \alg{T}^{i}
			}{
				a_{i} \in R,\
				a_{0} \in R^{\times}
			}.
	\]
Likewise, the multiplicative group $K^{\times}$ of $K$
can be viewed as a perfect group scheme $\alg{K}^{\times}$,
which is the direct product of $\alg{U}_{K}$
and the discrete infinite cyclic group generated by $\alg{T}$.
Let $K_{p}$ be the perfect closure of $K$.
We consider the $K_{p}$-rational point $1 - T \alg{T}^{-1}$ of $\alg{K}^{\times}$
and the corresponding morphism
	\[
		\varphi \colon \Spec K_{p} \to \alg{K}^{\times}.
	\]
This morphism induces a homomorphism
	\[
			\eta
		\colon
			\Gal(K^{\ab} / K)
		=
			\pi_{1}^{\et}(\Spec K_{p})^{\ab}
		\overset{\varphi}{\to}
			\pi_{1}^{k}(\alg{K}^{\times})
	\]
on the fundamental groups.
Here $K^{\ab}$ is the maximal abelian extension of $K$
and $\pi_{1}^{\et}(\cdot)^{\ab}$ denotes the maximal abelian quotient of the \'etale fundamental group.
The group $\pi_{1}^{k}(\alg{K}^{\times})$
is the fundamental group of $\alg{K}^{\times}$ as a perfect group scheme over $k$,
which classifies all surjective isogenies
to $\alg{K}^{\times}$ with finite constant kernels.
Now we state the main theorem of this paper:

\begin{mainth}
	\label{Thm:Main} \mbox{}
	\begin{enumerate}
		\item
			\label{Ass:Main:Isom}
			The above defined map
			$\eta \colon \Gal(K^{\ab} / K) \to \pi_{1}^{k}(\alg{K}^{\times})$
			is an isomorphism.
		\item
			\label{Ass:Main:Inverse}
			The inverse of $\eta$ restricted to $\pi_{1}^{k}(\alg{U}_{K})$
			coincides with the reciprocity isomorphism
				$
						\theta
					\colon
						\pi_{1}^{k}(\alg{U}_{K})
					\isomto
						I(K^{\ab} / K)
				$
			of Serre and Hazewinkel (\cite{Ser61}, \cite{DG70}),
			where $I$ denotes the inertia group.
	\end{enumerate}
\end{mainth}

Note that Assertion \ref{Ass:Main:Isom} of the theorem says that
we have an essentially one-to-one correspondence
between surjective isogenies $A \onto \alg{K}^{\times}$ with finite constant kernels
and finite abelian extensions $L$ of $K$ by pullback by $\varphi$,
as expressed as a cartesian diagram of the form
	\[
		\begin{CD}
				\Spec L_{p}
			@>>>
				A
			\\
			@VVV
			@VVV
			\\
				\Spec K_{p}
			@> \varphi >>
				\alg{K}^{\times}.
		\end{CD}
	\]
In other words, we have
$\Ext_{k}^{1}(\alg{K}^{\times}, N) \isomto H^{1}(K, N)$
for any finite constant group $N$.
Another remark is that
the proof of Assertion \ref{Ass:Main:Inverse} that we give in this paper
does not use the fact that $\theta$ is an isomorphism.
This means that Theorem \ref{Thm:Main} and our proof of it
together give another proof of this fact,
more specifically, the existence theorem of the local class field theory of Serre and Hazewinkel
(\cite[\S 4]{Ser61}, \cite[6.3]{DG70}).
For a generalization of Assertion \ref{Ass:Main:Inverse}
for the full groups $\pi_{1}^{k}(\alg{K}^{\times})$ and $\Gal(K^{\ab} / K)$,
see Remark \ref{Rmk:NotTR}.

We outline the paper.
After providing some preliminaries at Section \ref{Sect:Prelim},
we give three different proofs of Theorem \ref{Thm:Main}.
Each of them has its own advantages and interesting points.
The first proof given in Section \ref{Sect:DirectProof} may be the standard one.
The method is purely local.
In this proof, for Assertion \ref{Ass:Main:Isom},
we explicitly calculate the groups
$\Gal(K^{\ab} / K)$ and $\pi_{1}^{k}(\alg{K}^{\times})$
and the homomorphism $\eta$ between them individually.
For Assertion \ref{Ass:Main:Inverse},
we interpret the assertion to the compatibility of $\varphi$ with norm maps
and prove it by analyzing the diagonal divisors.
The second proof given in Section \ref{Sect:LTproof} relies on Lubin-Tate theory.
Hence it is applicable only for finite residue field cases.
The third proof given in Section \ref{Sect:GeomProof} is a geometric proof,
which was suggested by the referee to the author.
In this geometric proof, for Assertion \ref{Ass:Main:Isom},
we need the Albanese property of
the morphism $\varphi \colon \Spec K_{p} \to \alg{K}^{\times}$,
which was established by Contou-Carrere in \cite{CC81}
(see also \cite{CC94}).
For Assertion \ref{Ass:Main:Inverse},
we need the local-global compatibility
and the global version of Assertion \ref{Ass:Main:Inverse},
which was obtained by Serre in \cite[\S 5]{Ser61}.

In Section \ref{Sect:TwoDim},
we define a morphism analogous to
$\varphi \colon \Spec K_{p} \to \alg{K}^{\times}$
for a field of the form $k((S))((T))$ using the second algebraic $K$-group
instead of the multiplicative group.
If $k$ is a finite field,
the field $k((S))((T))$ is called
a two-dimensional local field (\cite{FK00}) of positive characteristic.
In analogy with the second proof of Theorem \ref{Thm:Main} by Lubin-Tate theory,
we attempt to formulate an analogous theory to Lubin-Tate theory
for the two-dimensional local field $k((S))((T))$.
Although we have not yet obtained an analog of a Lubin-Tate formal group in this paper,
we did obtain a result that may be thought of as a partial result
for abelian extensions having $p$-torsion Galois groups (Proposition \ref{Prop:TwoDim}).
In Section \ref{Sect:D-mod},
we give an analog of the result of \cite[\S 2.6]{BBDE} on $\D$-modules.
Their result is for fields of the form $k((T))$ with $k$ characteristic zero
while ours $k((S))((T))$ with $k$ characteristic zero.
This is also regarded as an analog of Proposition \ref{Prop:TwoDim}.

We give a couple of comments on literature.
The morphism $\varphi \colon \Spec K_{p} \to \alg{K}^{\times}$
had been introduced by Grothendieck (\cite[August 9, 1960]{CS01})
long before our paper.
Besides Contou-Carrere as mentioned above,
the morphism $\varphi$ has been studied by many people.
The $\D$-module version of Theorem \ref{Thm:Main} in the zero-characteristic case
had been established in \cite[\S 2.6]{BBDE}.
Deligne had given a sketch of proof of results
stronger than Assertion \ref{Ass:Main:Isom}.
This is written in Section e of
his letter to Serre contained in \cite{BE01}.
His method is different from our method.
The author did not realize these works
at the time of writing this paper.

\begin{enonce*}[remark]{Acknowledgements}
	This is an extended version of the master thesis of the author at Kyoto University.
	The author would like to express his deep gratitude
	to his advisor Kazuya Kato and Tetsushi Ito
	for their suggestion of the problems, encouragement, and many helpful discussions.
	The author is grateful also to the referee
	for the suggestion of the geometric proof of Theorem \ref{Thm:Main},
	and to Professor Takuro Mochizuki for informing the author of the work \cite{BBDE}.
\end{enonce*}


\section{Preliminaries}
\label{Sect:Prelim}

In this section, we give precise definitions of the terms
that we used in Introduction and will use in the following sections, fix notation,
prove the independence of the choice of the prime element $T$ (Proposition \ref{Prop:Indep}),
and reduce the proof of Theorem \ref{Thm:Main}
to the case of algebraically closed residue fields.


\subsection{Definitions and notation}
\label{Sect:DefNot}

We work on the site $\pfpqc{k}$
of perfect schemes over a perfect field $k$ of characteristic $p > 0$ with the fpqc topology
(compare with \cite[Chapter III, \S 0, ``Duality for unipotent perfect group schemes'']{Mil06};
see also \cite{SY12}).
The category of sheaves of abelian groups on $\pfpqc{k}$ contains
the category of commutative affine pro(-quasi)-algebraic groups over $k$
in the sense of Serre (\cite{Ser60})
as an thick abelian full subcategory.
We denote by $\Ext_{k}^{i}$ the $i$-th Ext functor
for the category of sheaves of abelian groups on $\pfpqc{k}$.
For a sheaf of abelian groups $A$ on $\pfpqc{k}$ and a non-negative integer $i$,
we define the fundamental group $\pi_{1}^{k}(A)$ of $A$
to be the Pontryagin dual of the torsion abelian group
$\injlim_{n \ge 1} \Ext_{k}^{1}(A, n^{-1} \Z / \Z)$.
If $A$ is an extension of an \'etale group
whose group of geometric points is finitely generated as an abelian group
by an affine proalgebraic group,
then $\injlim_{n \ge 1} \Ext_{k}^{1}(A, n^{-1} \Z / \Z) = \Ext_{k}^{1}(A, \Q / \Z)$
and hence the Pontryagin dual of $\Ext_{k}^{1}(A, \Q / \Z)$ coincides with $\pi_{1}^{k}(A)$.
For a $k$-algebra $R$, we denote by $R_{p}$
the perfect $k$-algebra given by
the injective limit of $p$-th power maps $R \to R \to \cdots$,
where the $i$-th copy of $R$ in this system is given the map
$k \to R$, $a \mapsto a^{p^{i}}$ as its $k$-algebra structure map.

Let $K$ be a complete discrete valuation field
of equal characteristic with residue field $k$.
We denote by $\Order_{K}$ the ring of integers of $K$
and by $\maxid_{K}$ the maximal ideal of $\Order_{K}$.
We set $U_{K} = U_{K}^{0} = \Order_{K}^{\times}$
and $U_{K}^{n} = 1 + \maxid_{K}^{n}$ for $n \ge 1$.
We fix an algebraic closure $\algcl{K}$ of $K$.
All algebraic extensions of $K$ are taken inside $\algcl{K}$.
We denote by $K^{\sep}$ ($\subset \algcl{K}$) the separable closure of $K$,
by $K^{\ur}$ the maximal unramified extension of $K$
and by $K^{\ab}$ the maximal abelian extension of $K$, respectively.
Since $K$ has equal characteristic $p > 0$,
the rings $\Order_{K}$ and $K$ have canonical structures of $k$-algebras
by the Teichm\"uller section $k \into \Order_{K}$.
Hence, for a sheaf of abelian groups $A$ on $\pfpqc{k}$,
we can consider the groups $A((\Order_{K})_{p})$ (resp.\ $A(K_{p})$)
of $(\Order_{K})_{p}$-rational (resp.\ $K_{p}$-rational) points of $A$.
We denote by $A((\maxid_{K})_{p})$ the kernel of the natural map
$A((\Order_{K})_{p}) \to A(k)$.
The reduction map $\Order_{K} \onto k$ and the Teichm\"uller section $k \into \Order_{K}$
give a canonical splitting $A((\Order_{K})_{p}) = A(k) \dsum A((\maxid_{K})_{p})$.

We define a sheaf of rings $\alg{O}_{K}$ on $\pfpqc{k}$
by setting $\alg{O}_{K}(R) = R \ctensor_{k} \Order_{K}$
for each perfect $k$-algebra $R$,%
\footnote{This defines a sheaf on the site of perfect $k$-algebras with the fpqc topology.
We actually need to take its Zariski sheafification to have a sheaf on $\pfpqc{k}$.
The two sites here have equivalent categories of sheaves.
The same process is applied to all sheaves in this paper.}
where $\ctensor$ denotes the completed tensor product.
Let $\alg{K}$ be the sheaf of rings on $\pfpqc{k}$
with $\alg{K}(R) = \alg{O}_{K}(R) \tensor_{\Order_{K}} K$.
We set $\alg{U}_{K} = \alg{O}_{K}^{\times}$.
For each $n \ge 0$, the sheaf of rings $\alg{O}_{K}$ has
a subsheaf of ideals $\alg{p}_{K}^{n}$ with
	$
			\alg{p}_{K}^{n}(R)
		=
			\alg{O}_{K}(R) \tensor_{\Order_{K}} \maxid_{K}^{n}
	$.
The presentation
	$
			\alg{O}_{K}
		=
			\projlim_{n \to \infty}
				\alg{O}_{K} / \alg{p}_{K}^{n}
	$
gives an affine proalgebraic ring structure for $\alg{O}_{K}$.
Likewise, $\alg{U}_{K}$ has a subsheaf of groups
$\alg{U}_{K}^{n} = 1 + \alg{p}_{K}^{n}$ for each $n \ge 1$
(for $n = 0$, we set $\alg{U}_{K}^{0} = \alg{U}_{K}$).
The presentation
	$
			\alg{U}_{K}
		=
			\projlim_{n \to \infty}
				\alg{U}_{K} / \alg{U}_{K}^{n}
	$
gives an affine proalgebraic group structure for $\alg{U}_{K}$.
We have a split exact sequence
$0 \to \alg{U}_{K} \to \alg{K}^{\times} \to \Z \to 0$.
For a prime element $T$ of $\Order_{K}$,
let $\alg{T}$ be the $k$-rational point of $\alg{O}_{K}$ defined by
$\alg{T} = 1 \tensor T \in k \ctensor_{k} \Order_{K} = \alg{O}_{K}(k)$.
If we fix a prime element $T$,
we may write $\alg{U}_{K}(R) = R[[\alg{T}]]^{\times}$ and
$\alg{K}^{\times}(R) = R[[\alg{T}]][\alg{T}^{-1}]^{\times}$,
which are the descriptions of $\alg{U}_{K}$ and $\alg{K}^{\times}$ given in Introduction.
The rational point $1 - T \alg{T}^{-1} \in \alg{K}^{\times}(K_{p})$
taken in Introduction should therefore be understood to be
	$
			1 - (T \tensor 1) \tensor T^{-1}
		\in
			(
					(K_{p} \ctensor_{k} \Order_{K})
				\tensor_{\Order_{K}}
					K
			)^{\times}
		=
			\alg{K}^{\times}(K_{p})
	$.

For a rational point $f \in \alg{K}^{\times}(K_{p})$,
we denote by $\varphi_{f}$ the corresponding morphism of $k$-schemes
$\Spec K_{p} \to \alg{K}^{\times}$.
An extension $0 \to N \to A \to \alg{K}^{\times} \to 0$ with $N$ finite constant
pulls back by $\varphi_{f} \colon \Spec K_{p} \to \alg{K}^{\times}$ to an $N$-torsor on $\Spec K_{p}$.
This defines a homomorphism
$\Ext_{k}^{1}(\alg{K}^{\times}, N) \to H^{1}(K, N)$.
By Pontryagin duality, we have a homomorphism
$\Gal(K^{\ab} / K) \to \pi_{1}^{k}(\alg{K}^{\times})$,
which we denote by $\eta_{f}$.


\subsection{A subset of $\alg{K}^{\times}(K_{p})$}
In this subsection,
we define a subset $\AJ{K}$ of $\alg{K}^{\times}(K_{p})$
and prove some properties of it.
This subset contains $1 - T \alg{T}^{-1}$ for all prime elements $T$
and is convenient for proving both the independence of the prime element
(Proposition \ref{Prop:Indep})
and the compatibility with norm maps (Proposition \ref{Prop:BC}).

\begin{defi}
	\label{Def:AJ}
	We define $\AJ{K}$ to be the set of elements $f \in \alg{K}^{\times}(K_{p})$
	satisfying the following conditions:
	\begin{enumerate}
		\item \label{Cond:Ideal}
			$f$ is in $\alg{K}((\Order_{K})_{p})$
			and generates as an ideal of $\alg{K}((\Order_{K})_{p})$
			the kernel of the multiplication map
				$
						\alg{K}((\Order_{K})_{p})
					=
							((\Order_{K})_{p} \ctensor_{k} \Order_{K})
						\tensor_{\Order_{K}}
							K
					\to
						K_{p}
				$ and
		\item \label{Cond:Conv}
			the natural map $\alg{K}((\Order_{K})_{p}) \to \alg{K}(k)$ sends $f$ to $1$.
	\end{enumerate}
\end{defi}

\begin{prop} \mbox{}
	\label{Prop:AJmaps}
	\begin{enumerate}
		\item \label{Ass:AJmaps:Containment}
			The set $\AJ{K}$ contains $1 - T \alg{T}^{-1}$ for any prime $T$.
		\item \label{Ass:AJmaps:Action}
			The group
				$
						\alg{U}_{K}((\maxid_{K})_{p})
					=
						\Ker \bigl(
							\alg{U}_{K}((\Order_{K})_{p}) \to \alg{U}_{K}(k)
						\bigl)
				$
			acts by multiplication on $\AJ{K}$ transitively without fixed points.
		\item \label{Ass:AJmaps:Valuation}
			The valuation map $\alg{K}^{\times} \to \Z$ sends
			every element of $\AJ{K}$ to $-1$.
		\item \label{Ass:AJmaps:Diagram}
			For every $f \in \AJ{K}$, we have the following commutative diagram
			with exact rows:
				\[
					\begin{CD}
							0
						@>>>
							I(K^{\ab} / K)
						@>>>
							\Gal(K^{\ab} / K)
						@>>>
							\Gal(k^{\ab} / k)
						@>>>
							0
						\\
						@.
						@VV \eta_{f} V
						@VV \eta_{f} V
						@VV \wr V
						@.
						\\
							0
						@>>>
							\pi_{1}^{k}(\alg{U}_{K})
						@>>>
							\pi_{1}^{k}(\alg{K}^{\times})
						@>>>
							\pi_{1}^{k}(\Z)
						@>>>
							0,
					\end{CD}
				\]
			where the right vertical isomorphism
			$\Gal(k^{\ab} / k) \isomto \pi_{1}^{k}(\Z)$
			is the Pontryagin dual map of the natural isomorphism
			$H^{1}(k, \Q / \Z) \cong \Ext_{k}^{1}(\Z, \Q / \Z)$ times $-1$.
		\item \label{Ass:AJmaps:Ram}
			For each prime element $T$, the restricted homomorphism
			$\eta_{f} \colon I(K^{\ab} / K) \to \pi_{1}^{k}(\alg{U}_{K})$
			for $f = 1 - T \alg{T}^{-1} \in \AJ{K}$ factors as the composite of
			the inclusion $I(K^{\ab} / K) \into \Gal(K^{\ab} / K)$
			and the homomorphism $\Gal(K^{\ab} / K) \to \pi_{1}^{k}(\alg{U}_{K})$
			induced by the morphism $\Spec K_{p} \to \alg{U}_{K}$
			corresponding to the rational point $- T + \alg{T} \in \alg{U}_{K}(K_{p})$.
	\end{enumerate}
\end{prop}

\begin{proof}
	\ref{Ass:AJmaps:Containment}.
	Obvious.
	
	\ref{Ass:AJmaps:Action}.
	If $f_{1}$ and $f_{2}$ are in $\AJ{K}$,
	we have $f_{2} = u f_{1}$ for some
		$
				u
			\in
				\Ker \bigl( \alg{K}^{\times}((\Order_{K})_{p}) \to \alg{K}^{\times}(k) \bigr)
			=
				\alg{K}^{\times}((\maxid_{K})_{p})
		$.
	Such $u$ is unique since $\alg{K}((\Order_{K})_{p})$ is a domain.
	We have $\alg{K}^{\times}((\maxid_{K})_{p}) = \alg{U}_{K}((\maxid_{K})_{p})$
	since $\Z((\maxid_{K})_{p}) = 0$.
	Thus we have Assertion \ref{Ass:AJmaps:Action}.
	
	\ref{Ass:AJmaps:Valuation}.
	This follows from Assertions \ref{Ass:AJmaps:Containment} and \ref{Ass:AJmaps:Action}.
	
	\ref{Ass:AJmaps:Diagram}.
	By Assertion \ref{Ass:AJmaps:Valuation}, we have a commutative diagram
		\[
			\begin{CD}
					\Spec K_{p}
				@>>>
					\Spec k
				\\
				@VV \varphi_{f} V
				@VVV
				\\
					\alg{K}^{\times}
				@>>>
					\Z,
			\end{CD}
		\]
	where the right vertical morphism corresponds to $-1 \in \Z(k)$.
	The maps induced on the fundamental groups give the required commutative diagram.
	The top row of the diagram in Assertion \ref{Ass:AJmaps:Diagram} is obviously exact.
	The bottom row of is exact as well
	since the sequence $0 \to \alg{U}_{K} \to \alg{K}^{\times} \to \Z \to 0$ is split exact.
	
	\ref{Ass:AJmaps:Ram}.
	The prime element $T$ gives a splitting $\alg{K}^{\times} \cong \alg{U}_{K} \times \Z$.
	The rational point $1 - T \alg{T}^{-1} \in \alg{K}^{\times}(K_{p})$
	corresponds to the rational point
	$(- T + \alg{T}, - 1) \in \alg{U}_{K}(K_{p}) \times \Z(K_{p})$.
	This implies Assertion \ref{Ass:AJmaps:Ram}.
\end{proof}


\subsection{Independence of the choice of the prime element}
\label{Sect:Indep}

We show that the morphisms $\Spec K_{p} \to \alg{K}^{\times}$
corresponding to the rational points
$1 - T \alg{T}^{-1} \in \alg{K}^{\times}(K_{p})$ for various primes $T$
actually induce the same homomorphism
$\Gal(K^{\ab} / K) \to \pi_{1}^{k}(\alg{K}^{\times})$ on the fundamental groups.
More generally, we prove the following.

\begin{prop}
	\label{Prop:Indep}
	The homomorphism
	$\eta_{f} \colon \Gal(K^{\ab} / K) \to \pi_{1}^{k}(\alg{K}^{\times})$
	for $f \in \AJ{K}$ is independent of the choice of $f$.
\end{prop}

Thus we may write $\eta$ instead of $\eta_{f}$.
To prove the proposition,
we need the following lemma,
which will also be used in the proof of Proposition \ref{Prop:BC}
and in Section \ref{Sect:GeomProof}.

\begin{lemm}
	\label{Lem:SurjOnPrimes}
	Let $0 \to N \to A \to A' \to 0$ be an exact sequence
	of sheaves of abelian groups on $\pfpqc{k}$ with $N$ representable by an affine scheme.
	Then the induced map $A((\maxid_{K})_{p}) \to A'((\maxid_{K})_{p})$ is surjective.
\end{lemm}

\begin{proof}
	Taking cohomology, we have a commutative diagram with exact rows
		\[
			\begin{CD}
					0
				@>>>
					N((\Order_{K})_{p})
				@>>>
					A((\Order_{K})_{p})
				@>>>
					A'((\Order_{K})_{p})
				@>>>
					H_{\fpqc}^{1}(\Order_{K}, N)
				\\
				@. @VVV @VVV @VVV @VVV
				\\
					0
				@>>>
					N(k)
				@>>>
					A(k)
				@>>>
					A'(k)
				@>>>
					H_{\fpqc}^{1}(k, N),
			\end{CD}
		\]
	where $H_{\fpqc}$ denotes the cohomology with respect to the site $\pfpqc{k}$.
	
	We show that the last vertical homomorphism is an isomorphism.
	Since $N$ is affine, we can write $N = \projlim N_{\lambda}$,
	where each $N_{\lambda}$ is an affine algebraic group scheme over $k$
	viewed as a sheaf on $\pfpqc{k}$ (\cite[III, \S3, 7.5, Cor.]{DG70}).
	We have
		$
				H_{\fpqc}^{1}(\Order_{K}, N)
			\cong
				\projlim H_{\fpqc}^{1}(\Order_{K}, N_{\lambda})
		$.
	Since any infinitesimal group becomes zero over $\pfpqc{k}$,
	we can take $N_{\lambda}$ to be smooth.
	Then we have
		$
				H_{\fpqc}^{1}(\Order_{K}, N_{\lambda})
			\cong
				H_{\et}^{1}(\Order_{K}, N_{\lambda})
			\cong
				H_{\et}^{1}(k, N_{\lambda})
		$
	by \cite[Chapter III, Theorem 3.9 and Remark 3.11]{Mil80}.
	Thus $H_{\fpqc}^{1}(\Order_{K}, N) \cong H_{\fpqc}^{1}(k, N)$.
	
	On the other hand, the horizontal homomorphisms have natural sections
	corresponding to the Teichm\"uller section $k \into \Order_{K}$.
	These facts imply the required surjectivity.
\end{proof}

\begin{proof}[Proof of Proposition \ref{Prop:Indep}]
	By Assertion \ref{Ass:AJmaps:Action} of Proposition \ref{Prop:AJmaps},
	it is enough to show that
	$\eta_{f} \colon \Gal(K^{\ab} / K) \to \pi_{1}^{k}(\alg{K}^{\times})$ is a zero map
	for $f \in \alg{U}_{K}((\maxid_{K})_{p}) = \alg{K}^{\times}((\maxid_{K})_{p})$.
	Let $0 \to N \to A \to \alg{K}^{\times} \to 0$ be an exact sequence
	of sheaves of abelian groups on $\pfpqc{k}$ with $N$ finite constant.
	Then the morphism
	$\varphi_{f} \colon \Spec K_{p} \to \alg{K}^{\times}$
	lifts to a morphism $\Spec K_{p} \to A$
	by Lemma \ref{Lem:SurjOnPrimes}.
	This implies that $\eta_{f} = 0$.
\end{proof}


\subsection{Reduction to the case of algebraically closed residue fields}
\label{Sect:RedToAlgCl}

We assume that
Theorem \ref{Thm:Main} is true for algebraically closed $k$
and show that Theorem \ref{Thm:Main} is true for general perfect $k$.
By Assertion \ref{Ass:AJmaps:Diagram} of Proposition \ref{Prop:AJmaps},
the problem is reduced to that for the restricted part
$\eta \colon I(K^{\ab} / K) \to \pi_{1}^{k}(\alg{U}_{K})$.
For this, it is enough to see that
the $\Gal(\algcl{k} / k)$-coinvariants of $\Gal((K^{\ur})^{\ab} / K^{\ur})$
(for the conjugation action) is $I (K^{\ab} / K)$
and the $\Gal(\algcl{k} / k)$-coinvariants of $\pi_{1}^{\algcl{k}}(\alg{U}_{K})$
is $\pi_{1}^{k}(\alg{U}_{K})$.
The first assertion follows from the fact that
the natural surjection $\Gal(K^{\sep} / K) \onto \Gal(\algcl{k} / k)$
admits a section (\cite[Chapter II, \S 4.3, Exercises]{Ser02}).
For the second assertion, we use the natural spectral sequence
	$
			H^{i}(
				k,
				\Ext_{\algcl{k}}^{j}(
					\alg{U}_{K}, \Q / \Z
				)
			)
		\Rightarrow
			\Ext_{k}^{i + j}(\alg{U}_{K}, \Q / \Z)
	$.
We have
	$
			\Hom_{\algcl{k}}(\alg{U}_{K}, \Q / \Z)
		=
			0
	$
by the connectedness of $\alg{U}_{K}$.
Thus
	$
			H^{0}(
				k,
				\Ext_{\algcl{k}}^{1}(
					\alg{U}_{K}, \Q / \Z
				)
			)
		=
			\Ext_{k}^{1}(\alg{U}_{K}, \Q / \Z)
	$.
Hence $\Gal(\algcl{k} / k)$-coinvariants of $\pi_{1}^{\algcl{k}}(\alg{U}_{K})$
is $\pi_{1}^{k}(\alg{U}_{K})$.
Thus we get the result.


\section{First proof of Theorem \ref{Thm:Main}}
\label{Sect:DirectProof}


\subsection{Proof of the part ``$\eta$ is an isomorphism''}
\label{Sect:Isom}

In this subsection, we assume that the residue field $k$ is algebraically closed
and fix a prime element $T$ of $\Order_{K}$
to prove that
	$
			\eta
		\colon
			I(K^{\ab} / K) = \Gal(K^{\ab} / K)
		\to
			\pi_{1}^{k}(\alg{K}^{\times}) = \pi_{1}^{k}(\alg{U}_{K})
	$
is an isomorphism.
By Assertion \ref{Ass:AJmaps:Ram} of Proposition \ref{Prop:AJmaps},
this homomorphism is induced by the morphism
$\varphi^{0} \colon \Spec K_{p} \to \alg{U}_{K}$
given by the rational point $- T + \alg{T} \in \alg{U}_{K}(K_{p})$.
Both groups $\Gal(K^{\ab} / K)$ and $\pi_{1}^{k}(\alg{U}_{K})$
are profinite abelian groups.
The Pontryagin dual of $\Gal(K^{\ab} / K)$ is $H^{1}(K, \Q / \Z)$.
The Pontryagin dual of $\pi_{1}^{k}(\alg{U}_{K})$ is
$\Ext_{k}^{1}(\alg{U}_{K}, \Q / \Z)$.
Therefore the problem is equivalent to showing that
the dual map
	$
			\eta^{\vee}
		\colon
			\Ext_{k}^{1}(\alg{U}_{K}, \Q / \Z)
		\to
			H^{1}(K, \Q / \Z)
	$
is an isomorphism.
Since $\alg{U}_{K} \cong \Gm \times \alg{U}_{K}^{1}$,
we have
	\[
			\Ext_{k}^{1}(\alg{U}_{K}, \Q / \Z)
		\cong
				\Ext_{k}^{1}(\Gm, \Q / \Z)
			\dsum
				\Ext_{k}^{1}(\alg{U}_{K}^{1}, \Q / \Z).
	\]


\subsubsection{The prime-to-p part}

We show that
	$
			\eta^{\vee}
		\colon
			\Ext_{k}^{1}(\alg{U}_{K}, \Q / \Z)
		\to
			H^{1}(K, \Q / \Z)
	$
induces an isomorphism on the prime-to-$p$ parts.
It is enough to show that
	$
			\eta^{\vee}
		\colon
			\Ext_{k}^{1}(\alg{U}_{K}, n^{-1} \Z / \Z)
		\to
			H^{1}(K, n^{-1} \Z / \Z)
	$
is an isomorphism for each integer $n$ prime to $p$.
Since $k$ is algebraically closed,
$n^{-1} \Z / \Z$ is isomorphic as a Galois module over $k$ to
the group $\mu_{n}$ of the $n$-th roots of unity.
The group $\Ext_{k}^{1}(\alg{U}_{K}^{1}, \mu_{n})$ is zero
since the $n$-th power map is an automorphism on $\alg{U}_{K}^{1}$
while $\mu_{n}$ is killed by $n$.
The group $\Ext_{k}^{1}(\Gm, \mu_{n})$ is a cyclic group of order $n$
generated by the extension class
$0 \to \mu_{n} \to \Gm \overset{n}{\to} \Gm \to 0$.
The group $H^{1}(K, \mu_{n})$ is a cyclic group of order $n$
generated by the Kummer character $\sigma \mapsto \sigma((- T)^{1 / n}) / (- T)^{1 / n}$.
The morphism $\varphi^{0} \colon \Spec K_{p} \to \alg{U}_{K}$
followed by the projection $\alg{U}_{K} \onto \Gm$
corresponds to the rational point $- T$.
These generators correspond each other via the homomorphism
	$
			\eta^{\vee}
		\colon
			\Ext_{k}^{1}(\alg{U}_{K}, \mu_{n})
		\to
			H^{1}(K, \mu_{n})
	$
since we have a cartesian diagram
	\[
		\begin{CD}
				\Spec K_{p}((-T)^{1 / n})
			@>>>
				\Gm
			\\
			@VVV
			@VV n V
			\\
				\Spec K_{p}
			@>>>
				\Gm.
		\end{CD}
	\]
Hence
	$
			\eta^{\vee}
		\colon
			\Ext_{k}^{1}(\alg{U}_{K}, \mu_{n})
		\to
			H^{1}(K, \mu_{n})
	$
is an isomorphism for any integer $n \ge 1$ prime to $p$.


\subsubsection{The $p$-primary part}

We show that the homomorphism on the $p$-primary parts
	$
			\eta^{\vee}
		\colon
			\Ext_{k}^{1}(\alg{U}_{K}, \Q_{p} / \Z_{p})
		\to
			H^{1}(K, \Q_{p} / \Z_{p})
	$
is an isomorphism as well.
We need the following lemmas.

\begin{lemm}
	\label{Lem:AHexp}
	The Artin-Hasse exponential map
		$
				\prod_{p \nmid n \ge 1} W
			\to
				\alg{U}_{K}^{1}
		$
	sending
		\[
				a
			=
				(a_{n})_{p \nmid n \ge 1}
			\in
				\prod_{p \nmid n \ge 1} W
			\quad \text{with} \quad
				a_{n}
			=
				(a_{n 0}, a_{n 1}, \dots)
			\in
				W
		\]
	to
		\[
				\prod_{p \nmid n \ge 1,\, m \ge 0}
					F(a_{n m } \alg{T}^{n p^{m}})
			\in
				\alg{U}_{K}^{1}
		\]
	is an isomorphism of proalgebraic groups.
	Here we denote by $W$ the additive group of Witt vectors
	and set
		$
				F(t)
			=
				\exp(
					- \sum_{e \ge 0}
						t^{p^{e}} / p^{e}
				)
			\in
				\Z_{p}[[t]]
		$.
\end{lemm}

\begin{proof}
	See \cite[Chapter V, \S 3, 16]{Ser88}.
\end{proof}

\begin{lemm} \label{Lem:p-div}
	Let $f \colon A \to B$ be a homomorphism between abelian groups $A$ and $B$.
	If both $A$ and $B$ are $p$-divisible and $p$-power torsion,
	and $f$ induces an isomorphism on the $p$-torsion parts,
	then $f$ itself is an isomorphism.
\end{lemm}

The group $\Ext_{k}^{1}(\Gm, \Q_{p} / \Z_{p})$ is zero
since the $p$-th power map is an automorphism on $\Gm$ as a perfect group scheme
while $\Q_{p} / \Z_{p}$ is a union of $p$-power torsion groups.
The group
	$
			\Ext_{k}^{1}(\alg{U}_{K}^{1}, \Q_{p} / \Z_{p})
		=
			\pi_{1}^{k}(\alg{U}_{K}^{1})^{\vee}
	$
is $p$-divisible
since $\alg{U}_{K}^{1} \cong \prod_{p \nmid n \ge 1} W$ by Lemma \ref{Lem:AHexp}
and $\pi_{1}^{k}(W)$ is $p$-torsion-free by \cite[\S 8.5, Prop.\ 5]{Ser60}.
The group $H^{1}(K, \Q_{p} / \Z_{p})$ is $p$-divisible
since the largest pro-$p$ quotient of $\Gal(K^{\sep} / K)$ is pro-$p$ free
(\cite[Chapter II, \S 2.2, Corollary 1]{Ser02}).
Thus, using Lemma \ref{Lem:p-div}, we are reduced to showing that
	$
			\eta^{\vee}
		\colon
			\Ext_{k}^{1}(\alg{U}_{K}^{1}, \Z / p \Z)
		\to
			H^{1}(K, \Z / p \Z)
	$
is an isomorphism.
We have
		\[
				\Ext_{k}^{1}(\alg{U}_{K}^{1}, \Z / p \Z)
			\cong
				\Dsum_{p \nmid n \ge 1}
					\Ext_{k}^{1}(W, \Z / p \Z)
			\cong
				\Dsum_{p \nmid n \ge 1}
					\Ext_{k}^{1}(\Ga, \Z / p \Z)
			\cong
				\Dsum_{p \nmid n \ge 1}
					k.
		\]
Here the last isomorphism comes from the isomorphism
$k \isomto \Ext_{k}^{1}(\Ga, \Z / p \Z)$
sending an element $a \in k^{\times}$ to the extension class given by
	\begin{equation} \label{Eq:ASisog}
			0
		\longrightarrow
			\Z / p \Z
		\longrightarrow
			\Ga
		\overset{a^{-1} \wp}{\longrightarrow}
			\Ga
		\longrightarrow
			0,
	\end{equation}
where $\wp(x) = x^{p} - x$ (\cite[\S 8.3, Prop.\ 3]{Ser60}).
On the other hand, the map defined by
	\begin{align*}
				\Dsum_{p \nmid n \ge 1}
					k T^{-n}
		&	\to
				H^{1}(K, \Z / p \Z),
		\\
				a T^{-n}
		&	\mapsto
				\bigl(
						\sigma
					\mapsto
							\sigma(\wp^{-1}(a T^{-n}))
						-
							\wp^{-1}(a T^{-n})
				\bigr)
	\end{align*}
is an isomorphism by Artin-Schreier theory.
To show that
	\[
			\eta^{\vee}
		\colon
			\Ext_{k}^{1}(\alg{U}_{K}, \Z / p \Z)
		\cong
			\Dsum_{p \nmid n \ge 1}
				k
		\to
			\Dsum_{p \nmid n \ge 1}
				k T^{-n}
		\cong
			H^{1}(K, \Z / p \Z)
	\]
is an isomorphism,
we need to calculate the following morphism:
	\[
			\Spec K_{p}
		\to
			\alg{U}_{K}^{1} / (\alg{U}_{K}^{1})^{p}
		\cong
			\prod_{p \nmid n \ge 1}
				W / p W
		\cong
			\prod_{p \nmid n \ge 1}
				\Ga.
	\]
The morphism $\Spec K_{p} \to \alg{U}_{K}^{1} / (\alg{U}_{K}^{1})^{p}$
corresponds to the $K_{p}$-rational point $1 - T^{-1} \alg{T}$
of $\alg{U}_{K}^{1} / (\alg{U}_{K}^{1})^{p}$.
The isomorphism
	$
			\prod_{p \nmid n \ge 1}
				\Ga
		\isomto
			\alg{U}_{K}^{1} / (\alg{U}_{K}^{1})^{p}
	$
sends each element $(a_{n})_{p \nmid n \ge 1}$ of the left-hand side
to $\prod_{p \nmid n \ge 1} F(a_{n} \alg{T}^{n})$ of the right-hand side.

\begin{prop} \mbox{}
	\begin{enumerate}
		\item
			\label{Ass:Dlog:Invmap}
			The inverse of the isomorphism 
				$
						\prod_{p \nmid n \ge 1}
							\Ga
					\isomto
						\alg{U}_{K}^{1} / (\alg{U}_{K}^{1})^{p}
				$
			is given by the map
				\[
						\alg{U}_{K}^{1} / (\alg{U}_{K}^{1})^{p}
					\overset{\dlog}{\to}
						\prod_{n \ge 1}
							\Ga \alg{T}^{n} \dlog \alg{T}
					\overset{\alpha}{\to}
						\prod_{p \nmid n \ge 1}
							\Ga,
				\]
			where
				$
						\alpha(
							\sum_{n \ge 1}
								b_{n} \alg{T}^{n} \dlog \alg{T}
						)
					=
						(- b_{n} / n)_{p \nmid n \ge 1}
				$.
		\item
			\label{Ass:Dlog:RatPt}
			The rational point $1 - T^{-1} \alg{T}$
			corresponds to $(1 / (n T^{n}))_{p \nmid n \ge 1}$
			via the isomorphism
				$
						\alg{U}_{K}^{1} / (\alg{U}_{K}^{1})^{p} (K_{p})
					\cong
						\prod_{p \nmid n \ge 1}
							\Ga(K_{p})
				$.
		\item
			\label{Ass:Dlog:Morphism}
			The map $\Spec K_{p} \to \prod_{p \nmid n \ge 1} \Ga$
			gives the $K_{p}$-rational point $(1 / (n T^{n}))_{p \nmid n \ge 1}$
			of $\prod_{p \nmid n \ge 1} \Ga$.
	\end{enumerate}
\end{prop}

\begin{proof}
	\ref{Ass:Dlog:Invmap}.
	Using the identity
		$
				\dlog F(t)
			=
				- \sum_{e \ge 0}
					t^{p^{e}} \dlog t
		$,
	we have
		\begin{align*}
					\dlog \left(
						 \prod_{p \nmid n \ge 1}
							F(a_{n} \alg{T}^{n})
					\right)
			&	=
					- \sum_{
						\substack{
							e \ge 0 \\
							p \nmid n \ge 1
						}
					}
						(a_{n} \alg{T}^{n})^{p^{e}}
						\dlog(a_{n} \alg{T}^{n})
			\\
			&	=
					\sum_{
						\substack{
							e \ge 0 \\
							p \nmid n \ge 1
						}
					}
						(- n) (a_{n} \alg{T}^{n})^{p^{e}}
						\dlog \alg{T}.
		\end{align*}
	Thus the map $\alpha \circ \dlog$ sends
	$\prod_{p \nmid n \ge 1} F(a_{n} \alg{T}^{n})$
	to $(a_{n})_{p \nmid n \ge 1}$, as desired.
	
	\ref{Ass:Dlog:RatPt}.
	A simple calculation shows that
		$
				(\alpha \circ \dlog) (1 - T^{-1} \alg{T})
			=
				(1 / (n T^{n}))_{p \nmid n \ge 1}
		$.
	
	\ref{Ass:Dlog:Morphism}.
	This follows from Assertion \ref{Ass:Dlog:RatPt}.
\end{proof}

Now we calculate $\eta^{\vee}$.
Let $n \ge 1$ be an integer prime to $p$
and $a \ne 0$ be an element of $k$
regarded as an element of $\Dsum_{p \nmid n \ge 1} k$
by the inclusion $k \into \Dsum_{p \nmid n \ge 1} k$ into the $n$-th summand.
The corresponding extension of $\Ga$ is given by \eqref{Eq:ASisog}.
We have a cartesian diagram
	\[
		\begin{CD}
				\Spec K_{p}(\wp^{-1}(a / n T^{n}))
			@>>>
				\Ga
			\\
			@VVV
			@VV a^{-1} \wp V
			\\
				\Spec K_{p}
			@>>>
				\Ga.
		\end{CD}
	\]
Thus
	$
			\eta^{\vee}
		\colon
			\Dsum_{p \nmid n \ge 1}
				k
		\to
			\Dsum_{p \nmid n \ge 1}
				k T^{-n}
	$
preserves the direct factors and
the map on the $n$-th factor is given by multiplication by $1 / n$.
This shows that
	$
			\eta^{\vee}
		\colon
			\Ext_{k}^{1}(\alg{U}_{K}, \Z / p \Z)
		\to
			H^{1}(K, \Z / p \Z)
	$
is an isomorphism.

This completes the proof of Assertion \ref{Ass:Main:Isom} of Theorem \ref{Thm:Main}.


\subsection{Proof of the part ``$\eta^{-1} = \theta$''}

In this section, we assume that
the residue field $k$ is a general perfect field of characteristic $p > 0$
(see Remark \ref{Rmk:NotTR} below for the point this generality helps for).
To prove Assertion \ref{Ass:Main:Inverse} of Theorem \ref{Thm:Main},
we first quickly recall the construction of the reciprocity map
$\theta \colon \pi_{1}^{k}(\alg{U}_{K}) \isomto I(K^{\ab} / K)$
of the local class field theory of Serre and Hazewinkel.
Let $L$ be a finite totally ramified Galois extension of $K$.
We can apply the constructions in Section \ref{Sect:DefNot} to $L$
to get sheaves $\alg{U}_{L}$ and $\alg{L}^{\times}$ on $\pfpqc{k}$.
Let $\alg{V}_{L / K}$ be the subgroup of $\alg{U}_{L}$ generated by
elements of the form $\sigma(u) / u$ for $\sigma \in \Gal(L / K)$ and $u \in \alg{U}_{L}$.
By \cite[Appendice, 4.2]{DG70}, we have exact sequences
	\begin{gather} \label{Eq:NormExSeqU}
			0
		\to
			\Gal(L / K)^{\ab}
		\to
			\alg{U}_{L} / \alg{V}_{L / K}
		\overset{N_{L / K}}{\to}
			\alg{U}_{K}^{\times}
		\to
			0,
		\\ \label{Eq:NormExSeqM}
			0
		\to
			\Gal(L / K)^{\ab}
		\to
			\alg{L}^{\times} / \alg{V}_{L / K}
		\overset{N_{L / K}}{\to}
			\alg{K}^{\times}
		\to
			0,
	\end{gather}
where the second map $N_{L / K}$ is the norm map
and the first map sends $\sigma$ to $\sigma(\boldsymbol{\pi}_{L}) / \boldsymbol{\pi}_{L}$
for any choice of a prime element $\pi_{L}$ of $\Order_{L}$,
with $\boldsymbol{\pi}_{L}$ defined as
$1 \tensor \pi_{L} \in k \ctensor_{k} \Order_{L} = \alg{O}_{L}(k)$
and $\sigma$ acting on $\alg{L}^{\times}$ as a morphism of sheaves.
To write the isomorphism
$\theta \colon \pi_{1}^{k}(\alg{U}_{K}) \isomto I(K^{\ab} / K)$,
note that the Pontryagin dual groups of $\pi_{1}^{k}(\alg{U}_{K})$, $I(K^{\ab} / K)$ are
$\Ext_{k}^{1}(\alg{U}_{K}, \Q / \Z)$, $H^{1}(I(K^{\ab} / K), \Q / \Z)$, respectively.
Now the dual map
	$
			\theta^{\vee}
		\colon
			H^{1}(I(K^{\ab} / K), G)
		\isomto
			\Ext_{k}^{1}(\alg{U}_{K}, G)
	$
for a finite abelian group $G$
is given by sending $L$ totally ramified abelian with $\Gal(L / K) = G$
to the extension class \eqref{Eq:NormExSeqU}
(\cite[Appendice, 7.1]{DG70}).

On the other hand, the dual
	$
			\eta^{\vee}
		\colon
			\Ext_{k}^{1}(\alg{K}^{\times}, G)
		\isomto
			H^{1}(K, G)
	$
of our isomorphism
$\eta \colon \Gal(K^{\ab} / K) \isomto \pi_{1}^{k}(\alg{K}^{\times})$
sends the extension class \eqref{Eq:NormExSeqM} to
the $G$-torsor on $\Spec K_{p}$ given by pulling back \eqref{Eq:NormExSeqM} by the morphism
$\varphi_{f} \colon \Spec K_{p} \to \alg{K}^{\times}$, where $f \in \AJ{K}$.
Therefore our claim $\eta^{-1} = \theta$ is interpreted as the claim that
the $G$-torsor on $\Spec K_{p}$ thus obtained is $\Spec L_{p}$.
Hence Assertion \ref{Ass:BCwActions} of the following proposition proves
Assertion \ref{Ass:Main:Inverse} of Theorem \ref{Thm:Main}:

\begin{prop} \label{Prop:BC}
	Assume that the residue field $k$ is perfect of characteristic $p > 0$.
	Let $L / K$ be a finite totally ramified Galois extension.
	\begin{enumerate}
	\item \label{Ass:BCUptoUnits}
		For any $f \in \AJ{K}$ and any $g \in \AJ{K}$,
		the diagram
			\[
				\begin{CD}
						\Spec L_{p}
					@>> \varphi_{g} >
						\alg{L}^{\times}
					\\
					@VVV
					@VV N_{L / K} V
					\\
						\Spec K_{p}
					@> \varphi_{f} >>
						\alg{K}^{\times}
				\end{CD}
			\]
		is commutative ``up to a $(\maxid_{L})_{p}$-rational point'',
		namely the ratio of the two elements of $\alg{K}^{\times}(L_{p})$
		corresponding to the two morphisms
		$\Spec L_{p} \rightrightarrows \alg{K}^{\times}$
		coming from the above possibly non-commutative diagram
		is in
			$
					\alg{K}^{\times}((\maxid_{L})_{p})
				=
					\Ker(
							\alg{K}^{\times}(\Order_{L})_{p})
						\onto
							\alg{K}^{\times}(k)
					)
			$.
	\item \label{Ass:BCComm}
		For any $f \in \AJ{K}$,
		there exists $g \in \AJ{L}$ such that the above diagram is indeed commutative.
	\item \label{Ass:BCwActions}
		Assume that $L / K$ is abelian.
		For any $f \in \AJ{K}$ and any $g \in \AJ{L}$ that makes the diagram commutative,
		the induced diagram
			\[
				\begin{CD}
						\Spec L_{p}
					@>> \varphi_{g} >
						\alg{L}^{\times} / \alg{V}_{L / K}
					\\
					@VVV
					@VV N_{L / K} V
					\\
						\Spec K_{p}
					@> \varphi_{f} >>
						\alg{K}^{\times}
				\end{CD}
			\]
		is cartesian, and is equivariant under the action of $\Gal(L / K)$ on the right vertical arrow
		coming from the exact sequence \eqref{Eq:NormExSeqM}
		and the natural action of $\Gal(L / K)$ on the left vertical arrow.
	\end{enumerate}
\end{prop}

\begin{proof}
	\ref{Ass:BCUptoUnits}.
	The two elements of $\alg{K}^{\times}(L_{p})$ mentioned in Assertion \ref{Ass:BCUptoUnits}
	are given by $f$ and $N_{L / K} g$.
	It is enough to show that $f / N_{L / K} g \in \alg{K}^{\times}((\Order_{L})_{p})$
	since then we immediately have $f / N_{L / K} g \in \alg{K}^{\times}((\maxid_{L})_{p})$
	by Condition \ref{Cond:Conv} of Definition \ref{Def:AJ}.
	The claim $f / N_{L / K} g \in \alg{K}^{\times}((\Order_{L})_{p})$ is equivalent to
	the equality $(N_{L / K} g) = (f)$ of ideals of the ring $\alg{K}((\Order_{L})_{p})$.
	Let $\Delta_{K} \subset \Spec \alg{K}((\Order_{K})_{p})$ be the diagonal divisor
	defined by the kernel of the multiplication map
		\[
				\alg{K}((\Order_{K})_{p})
			=
					((\Order_{K})_{p} \ctensor_{k} \Order_{K})
				\tensor_{\Order_{K}}
					K
			\to
				K_{p}.
		\]
	Then we have $\Delta_{K} = (f)$ as divisors on $\Spec \alg{K}((\Order_{K})_{p})$
	by Condition \ref{Cond:Ideal} of Definition \ref{Def:AJ}.
	Similarly, for $\sigma \in \Gal(L / K)$,
	let $\sigma(\Delta_{L}) \subset \Spec \alg{L}((\Order_{L})_{p})$ be
	the $\sigma$-twist of the diagonal divisor defined by the kernel of the map
		\[
				\alg{L}((\Order_{L})_{p})
			=
					((\Order_{L})_{p} \ctensor_{k} \Order_{L})
				\tensor_{\Order_{L}}
					L
			\overset{(1 \tensor \sigma) \tensor \sigma}{\longrightarrow}
				L_{p}.
		\]
	Then we have $p_{\ast} \Delta_{L} = (N_{L / K} g)$ as divisors on $\Spec \alg{K}((\Order_{L})_{p})$,
	where $p \colon \Spec \alg{L}((\Order_{L})_{p}) \to \Spec \alg{K}((\Order_{L})_{p})$ is the natural map.
	Therefore what we have to show is the equality $p_{\ast} \Delta_{L} = \Delta_{K}$
	of divisors on $\Spec \alg{K}((\Order_{L})_{p})$.
	This equality follows from the equalities
		\[
				p_{\ast} \Delta_{L}
			=
				\sum_{\sigma \in \Gal(L / K)} \sigma(\Delta_{L})
			=
				\Delta_{K}
		\]
	as divisors on $\Spec \alg{L}((\Order_{L})_{p})$.
	This proves Assertion \ref{Ass:BCUptoUnits}.
	
	\ref{Ass:BCComm}.
	First take $g \in \AJ{L}$ to be arbitrary.
	Then by Assertion \ref{Ass:BCUptoUnits}, we have
	$f / N_{L / K} g \in \alg{K}^{\times}((\maxid_{L})_{p})$.
	The homomorphism $N_{L / K} \colon \alg{L}^{\times} \to \alg{K}^{\times}$
	and its restriction $\alg{U}_{L} \to \alg{U}_{K}$ are surjective and have the same kernel;
	see \eqref{Eq:NormExSeqU} and \eqref{Eq:NormExSeqM}.
	Since $\alg{U}_{L}$ and $\alg{U}_{K}$ are affine, so is $\Ker(N_{L / K})$.
	Thus we can choose an element $u \in \alg{L}^{\times}((\maxid_{L})_{p})$
	so that $N_{L / K} u = f / N_{L / K} g$
	by Lemma \ref{Lem:SurjOnPrimes}.
	Replacing $g$ by $u g$, we have $N_{L / K} g = f$.
	This choice of $g$ makes the diagram commutative.
	
	\ref{Ass:BCwActions}.
	Before the proof, note that the field $\alg{L}(L_{p})$ has two different actions of $\Gal(L / K)$.
	One comes from the action on the sheaf $\alg{L}$,
	which we have already denoted by $h \mapsto \sigma(h)$, $\sigma \in \Gal(L / K)$.
	The other comes from the action on the coefficient field $L_{p}$,
	which we denote by $h \mapsto [\sigma](h)$, $\sigma \in \Gal(L / K)$.
	
	Now we show the equivariance first.
	This is the same, in terms of the rational point $g$, as the equality
	$[\sigma](g) = (\sigma(\boldsymbol{\pi}_{L}) / \boldsymbol{\pi}_{L}) \cdot g$
	in $(\alg{L}^{\times} / \alg{V}_{L / K})(L_{p})$.
	Below we use the exponential notation,
	so the equality to be shown can be written as
	$g^{[\sigma] - 1} = \boldsymbol{\pi}_{L}^{\sigma - 1}$.
	We have
		\[
				g^{[\sigma] - 1}
			=
				(g^{- [\sigma]})^{\sigma - 1}
				g^{\sigma [\sigma] - 1}.
		\]
	By Assertion \ref{Ass:AJmaps:Valuation} of Proposition \ref{Prop:AJmaps},
	the element $g^{- [\sigma]}$ is a prime element of the complete discrete valuation ring
	$\alg{O}_{L}(L_{p})$.
	Hence $(g^{- [\sigma]})^{\sigma - 1} = \boldsymbol{\pi}_{L}^{\sigma - 1}$
	in $(\alg{L}^{\times} / \alg{V}_{L / K})(L_{p})$.
	Therefore it is enough to show that
	$g^{\sigma [\sigma] - 1} \in \alg{V}_{L / K}(L_{p})$.
	
	By Condition \ref{Cond:Ideal} of Definition \ref{Def:AJ},
	the element $g^{\sigma [\sigma]}$ generates as an ideal of $\alg{L}((\Order_{L})_{p})$
	the kernel of the map
		\[
				\alg{L}((\Order_{L})_{p})
			=
					((\Order_{L})_{p} \ctensor_{k} \Order_{L})
				\tensor_{\Order_{L}}
					L
			\overset{(\sigma \tensor \sigma) \tensor \sigma}{\longrightarrow}
				L_{p},
		\]
	which is the same as the kernel of the usual multiplication map
		$
				\alg{L}((\Order_{L})_{p})
			=
					((\Order_{L})_{p} \ctensor_{k} \Order_{L})
				\tensor_{\Order_{L}}
					L
			\to
				L_{p}
		$.
	Therefore we have $g^{\sigma [\sigma] - 1} \in \alg{L}^{\times}((\Order_{L})_{p})$.
	Since $g \in \alg{L}^{\times}(L_{p})$,
	we have $g^{\sigma [\sigma] - 1} \in \alg{U}_{L}((\Order_{L})_{p})$.
	Moreover, we have $g^{\sigma [\sigma] - 1} \in \alg{U}_{L}((\maxid_{L})_{p})$
	by Condition \ref{Cond:Conv} of Definition \ref{Def:AJ}.
	Furthermore, we have $g^{\sigma [\sigma] - 1} \in (\Ker N_{L / K})((\maxid_{L})_{p})$
	since $f$ is $K$-rational and therefore
	$N_{L / K}(g^{\sigma [\sigma] - 1}) = f^{[\sigma] - 1} = 1$.
	The exact sequence \eqref{Eq:NormExSeqU} leads an exact sequence
	$0 \to \alg{V}_{L / K} \to \Ker N_{L / K} \to \Gal(L / K)^{\ab} \to 0$,
	hence $(\Ker N_{L / K})((\maxid_{L})_{p}) = \alg{V}_{L / K}((\maxid_{L})_{p})$.
	Thus we have $g^{\sigma [\sigma] - 1} \in \alg{V}_{L / K}((\maxid_{L})_{p})$.
	This proves the equivariance of the diagram.
	
	Finally we show that the diagram is cartesian.
	The pullback of the exact sequence \eqref{Eq:NormExSeqM} by
	$\varphi_{f} \colon \Spec K_{p} \to \alg{K}^{\times}$
	gives a $\Gal(L / K)$-torsor $Z$ on $\Spec K_{p}$.
	The diagram factors through a morphism from $\Spec L_{p}$ to $Z$ by universality.
	The equivariance just proved says that
	this is a morphism of $\Gal(L / K)$-torsors on $\Spec K_{p}$.
	Such a morphism is automatically an isomorphism.
	Therefore $Z = \Spec L_{p}$ and the diagram is cartesian.
\end{proof}

\begin{rema}
	\label{Rmk:NotTR}
	In \cite[\S 4]{SY10},
	the reciprocity isomorphism
	$\theta \colon I(K^{\ab} / K) \isomto \pi_{1}^{k}(\alg{U}_{K})$
	of Serre and Hazewinkel was extended to an isomorphism
	$\Gal(K^{\ab} / K) \isomto \pi_{1}^{k}(\alg{K}^{\times})$.
	A full discussion about this result was taken place in \cite{SY12}.
	The dual map
		$
				\theta^{\vee}
			\colon
				H^{1}(K, G)
			\isomto
				\Ext_{k}^{1}(\alg{K}^{\times}, G)
		$
	sends a finite totally ramified abelian extension $L$ of Galois group $G$
	to the extension class \eqref{Eq:NormExSeqM}.
	The above proposition also proves the equality $\eta^{-1} = \theta$
	as isomorphisms between full $\Gal(K^{\ab} / K)$ and $\pi_{1}^{k}(\alg{K}^{\times})$.
\end{rema}


\section{Second proof of Theorem \ref{Thm:Main}: Lubin-Tate theory}
\label{Sect:LTproof}

In this section, we assume that the residue field $k$ is finite
and give another proof of Theorem \ref{Thm:Main} for finite $k$ using Lubin-Tate theory.
Assume $k$ is the finite field with $q$ elements.
We fix a prime element $T$ of $\Order_{K}$.
This defines a splitting $\alg{K}^{\times} \cong \alg{U}_{K} \times \Z$.
By Assertion \ref{Ass:AJmaps:Ram} of Proposition \ref{Prop:AJmaps}, the homomorphism
	$
			\Gal(K^{\ab} / K)
		\overset{\eta}{\to}
			\alg{K}^{\times}
		\onto
			\alg{U}_{K}
	$
comes from the morphism $\varphi^{0} \colon \Spec K_{p} \to \alg{U}_{K}$
corresponding the rational point $- T + \alg{T} \in \alg{U}_{K}(K_{p})$.
We consider the short exact sequence
$0 \to U_{K} \to \alg{U}_{K} \stackrel{F - 1}{\to} \alg{U}_{K} \to 0$,
where $F$ is the $q$-th power Frobenius morphism.
The isogeny $F - 1$, called the Lang isogeny,
is universal among isogenies onto $\alg{U}_{K}$ with (pro-) finite constant kernels
(\cite[Chap.\ VI, \S 1, Prop.\ 6]{Ser88}),
hence it identifies $U_{K}$ as $\pi_{1}^{k}(\alg{U}_{K})$.

\begin{prop}
	\label{Prop:LT}
	There is a cartesian diagram
		\[
			\begin{CD}
					\Spec (K_{T}^{\ram})_{p}
				@>>>
					\alg{U}_{K}
				\\
				@VVV
				@VV F - 1 V
				\\
					\Spec K_{p}
				@> \varphi^{0} >>
					\alg{U}_{K}.
			\end{CD}
		\]
	Here $K_{T}^{\ram}$ is the field $K$
	adjoining all the $T^{m}$-torsion points
	(where $m$ runs through the integers $\ge 1$)
	of the Lubin-Tate formal group $F_{f}$ (\cite{Iwa86})
	whose equation of formal multiplication by $T$ is equal to $f(X) = T X + X^{q}$.
	The morphism $\Spec (K_{T}^{\ram})_{p} \to \alg{U}_{K}$ corresponds to
	the rational point
		$
			\sum_{m = 0}^{\infty}
				\alpha_{m + 1} \alg{T}^{m}
		$,
	where $\alpha_{m}$ is a generator of the module of $T^{m}$-torsion points of $F_{f}$
	such that $f(\alpha_{m + 1}) = \alpha_{m}$.
	The induced isomorphism
	$\Gal(K_{T}^{\ram} / K) \cong U_{K}$
	coincides with the one given by Lubin-Tate theory.
\end{prop}

\begin{proof}
	We calculate the geometric fiber of $F - 1$ over $- T + \alg{T}$.
	Let $g = \sum a_{m} \alg{T}^{m}$ be an element of $\alg{U}_{K}(\algcl{K})$.
	The equation $F(g) / g = - T + \alg{T}$
	is equivalent to the system of equations
		$f(a_{0}) = 0$,
		$f(a_{m + 1}) = a_{m}$,
		$m \ge 0$.
	Thus, for each $m \ge 0$, $a_{m}$ is a generator
	of the module of $T^{m + 1}$-torsion points of $F_{f}$.
	This proves the existence of the above cartesian diagram.
	Next we calculate the action of $\Gal(K_{T}^{\text{ram}} / K)$
	on the fiber of $F - 1$ over $- T + \alg{T}$.
	The Lubin-Tate group $F_{f}$ for $f(X) = T X + X^{q}$
	is the formal completion $\fGa$ of the additive group
	with formal multiplication of each element $\sum b_{m} T^{m}$ of $\Order_{K}$
	being given by the power series $\sum b_{m} f^{\circ m}(X) \in \End \fGa$,
	where $f^{\circ m}$ is the $m$-th iteration of $f$.
	Hence, if $\sigma$ corresponds to $u(T) = \sum b_{m} T^{m}$
	via the isomorphism
		$
				\Gal(K_{T}^{\text{ram}} / K)
			\cong
				U_{K}
		$
	of Lubin-Tate theory, we have
		\begin{align*}
				\sigma \left(
					\sum_{m \ge 0}
						\alpha_{m + 1} \alg{T}^{m}
				\right)
		&	=
				\sum_{m \ge 0}
					\sigma(\alpha_{m + 1}) \alg{T}^{m}
			=
				\sum_{0 \le k \le m < \infty}
					b_{k} \alpha_{m + 1 - k} \alg{T}^{m}
		\\
		&	=
				u(\alg{T})
				\sum_{m \ge 0}
					\alpha_{m + 1} \alg{T}^{m}.
		\end{align*}
	Thus the action of $\sigma$ on the fiber of $F - 1$ over $- T + \alg{T}$
	is given by multiplication by $u(\alg{T})$, as required.
\end{proof}

Therefore the map $\eta \colon \Gal(K^{\ab} / K) \to \pi_{1}^{k}(\alg{K}^{\times})$
factors through the restriction map
$\Gal(K^{\ab} / K) \onto \Gal(K_{T}^{\ram} K^{\ur} / K)$,
the isomorphism of Lubin-Tate theory
$\Gal(K_{T}^{\ram} K^{\ur} / K) \isomto (K^{\times})^{\wedge}$
(the profinite completion of $K^{\times}$)
and the natural isomorphism
$(K^{\times})^{\wedge} \cong \pi_{1}^{k}(\alg{K}^{\times})$.
Thus the assertion that $\eta$ is an isomorphism
is equivalent to the local Kronecker-Weber theorem for Lubin-Tate extensions:
$K^{\ab} = K_{T}^{\ram} K^{\ur}$.
Since the isomorphism $\theta$ of Serre and Hazewinkel for finite $k$
coincides with the one given by Lubin-Tate theory,
the equality $\eta^{-1} = \theta$ for such $k$ follows.

\begin{rema}
	This proof also shows that
	the proof for $\eta$ being an isomorphism
	given in Section \ref{Sect:Isom}
	gives another proof of the local Kronecker-Weber theorem for Lubin-Tate extensions
	in the equal characteristic case.
\end{rema}


\section{Third proof of Theorem \ref{Thm:Main}: geometry}
\label{Sect:GeomProof}

In this section, we give a geometric proof of Theorem \ref{Thm:Main}.
This method was suggested by the referee to the author.
As in Section \ref{Sect:DirectProof},
we assume the residue field $k$ is algebraically closed and
fix a prime element $T$ of $\Order_{K}$.
Both the morphism $\varphi \colon \Spec K_{p} \to \alg{K}^{\times}$,
which corresponds to the rational point $1 - T \alg{T}^{-1}$,
and the morphism $\varphi^{0} \colon \Spec K_{p} \to \alg{U}_{K}$,
which corresponds to the rational point $- T + \alg{T}$,
give the same homomorphism
	$
			\eta
		\colon
			\Gal(K^{\ab} / K)
		\to
			\pi_{1}^{k}(\alg{K}^{\times})
		=
			\pi_{1}^{k}(\alg{U}_{K})
	$.
We use the same letter to express morphisms and corresponding rational points.
It is enough to study the part
	$
			\eta^{\vee}
		\colon
			\Ext_{k}^{1}(\alg{U}_{K}, N)
		\to
			H^{1}(K, N)
	$
for each finite constant group $N$.

For the proof of Assertion \ref{Ass:Main:Isom} of Theorem \ref{Thm:Main},
we need the following result due to Contou-Carrere (\cite{CC81}).
Let $A$ be a commutative quasi-algebraic group over $k$.
Then the homomorphism
	\[
		\Hom(\alg{U}_{K}, A) \to A(K_{p}) / A((\Order_{K})_{p})
	\]
induced by the morphism $\varphi^{0} \colon \Spec K_{p} \to \alg{U}_{K}$
is an isomorphism (\cite[Th\'eor\`em (1.5)]{CC81}).
If $A$ is a commutative proalgebraic group over $k$,
the homomorphism $\Hom(\alg{U}_{K}, A) \to A(K_{p}) / A((\Order_{K})_{p})$
induced by the morphism $\varphi^{0} \colon \Spec K_{p} \to \alg{U}_{K}$
is still an isomorphism.

First we show that
	$
			\eta^{\vee}
		\colon
			\Ext_{k}^{1}(\alg{U}_{K}, N)
		\to
			H^{1}(K, N)
	$
is surjective, that is,
every finite abelian extension $L$ of $K$ with Galois group $N$
is obtained by the pullback of an extension $0 \to N \to A \to \alg{U}_{K} \to 0$
by the morphism $\varphi^{0} \colon \Spec K \to \alg{U}_{K}$.
Since $k$ is algebraically closed,
the abelian extension $L / K$ is
a composite of a Kummer extension of degree prime to $p$
and Artin-Schreier-Witt extensions.
A Kummer extension of $K$ of degree $n$, $p \nmid n$,
is the pullback of $\Gm \overset{n}{\to} \Gm$
by some morphism $\Spec K_{p} \to \Gm$.
An Artin-Schreier-Witt extension of $K$ of degree $p^{n}$
is the pullback of $W_{n} \overset{\wp}{\to} W_{n}$
by some morphism $\Spec K_{p} \to W_{n}$,
where $W_{n}$ is the additive group of Witt vectors of length $n$,
$\wp(x) = F x - x$ and $F$ is the Frobenius.
Combining them, we can find an exact sequence
$0 \to N \to A' \overset{\alpha}{\to} A'' \to 0$
of commutative quasi-algebraic groups
and a cartesian diagram
	\[
		\begin{CD}
				\Spec L_{p}
			@>> \beta >
				A'
			\\
			@VVV
			@V \alpha VV
			\\
				\Spec K_{p}
			@> \gamma >>
				A''.
		\end{CD}
	\]
Since $\varphi^{0} \colon \Spec K_{p} \to \alg{U}_{K}$ induces an isomorphism
$\Hom(\alg{U}_{K}, A'') \isomto A''(K_{p}) / A''((\Order_{K})_{p})$
by the result of Contou-Carrere quoted above,
we have $\gamma + \varepsilon = \varphi^{0} \circ \delta$
for some $\delta \colon \alg{U}_{K} \to A''$ and $\varepsilon \in A''((\Order_{K})_{p})$.
Since $\alpha \colon A' \to A''$ is surjective
and $k$ is algebraically closed,
Lemma \ref{Lem:SurjOnPrimes} shows that
$\alpha \colon A'((\Order_{K})_{p}) \to A''((\Order_{K})_{p})$ is also surjective.
Hence $\varepsilon = \alpha(\zeta)$ for some $\zeta \in A'((\Order_{K})_{p})$.
The above diagram is still commutative and cartesian
if we add $\zeta$ to $\beta$ and $\varepsilon$ to $\gamma$.
The resulting cartesian diagram splits into
the following commutative diagram with cartesian squares:
	\[
		\begin{CD}
				\Spec L_{p}
			@>>>
				A
			@>>>
				A'
			\\
			@VVV
			@VVV
			@V \alpha VV
			\\
				\Spec K_{p}
			@> \varphi^{0} >>
				\alg{U}_{K}
			@> \gamma >>
				A''
		\end{CD}
	\]
Thus we get the surjectivity of
	$
			\eta^{\vee}
		\colon
			\Ext_{k}^{1}(\alg{U}_{K}, N)
		\to
			H^{1}(K, N)
	$.

Next we show that $\eta^{\vee}$ is injective.
Let $0 \to N \to A \overset{\alpha}{\to} \alg{U}_{K} \to 0$ be an extension
whose image in $H^{1}(K, N)$ by $\eta^{\vee}$ is trivial,
i.e. $\alpha \circ \beta = \varphi^{0}$ for some $\beta \in A(K_{p})$.
The diagram
	\[
		\begin{CD}
				\Hom(\alg{U}_{K}, A)
			@=
				 A(K_{p}) / A((\Order_{K})_{p})
			\\
			@VV \alpha V
			@V \alpha VV
			\\
				\Hom(\alg{U}_{K}, \alg{U}_{K})
			@=
				\alg{U}_{K}(K_{p}) / \alg{U}_{K}((\Order_{K})_{p})
		\end{CD}
	\]
is commutative.
Let $\gamma \in \Hom(\alg{U}_{K}, A)$ correspond $\beta \in A(K_{p})$
via the top horizontal map.
Then we have $\alpha \circ \gamma = \id$
by the equality $\alpha \circ \beta = \varphi^{0}$ and this diagram.
Hence the extension $0 \to N \to A \overset{\alpha}{\to} \alg{U}_{K} \to 0$
is trivial.
Thus we get the injectivity of
	$
			\eta^{\vee}
		\colon
			\Ext_{k}^{1}(\alg{U}_{K}, N)
		\to
			H^{1}(K, N)
	$.

Finally we show that $\eta^{-1} = \theta$.
For this, we need generalized Jacobians (\cite{Ser88}).
We denote by $\Proj^{1}$ the projective line over $k$ with coordinate $T$.
For each effective divisor $\mathfrak{m}$ on $\Proj^{1}$
with support $\Supp{\mathfrak{m}} \subset \Proj^{1}$,
let $J_{\Proj^{1}, \mathfrak{m}}$ be
the moduli of line bundles on $X$ with level $\mathfrak{m}$ structure,
which is isomorphic to the direct product of $\Z$ and
the generalized Jacobian of $\Proj^{1}$ with modulus $\mathfrak{m}$.
Let
	$
			\Phi_{\mathfrak{m}}
		\colon
			\Proj^{1} - \Supp{\mathfrak{m}}
		\to
			J_{\Proj^{1}, \mathfrak{m}}
	$
be the Abel-Jacobi map times $-1$.
Since $\Proj^{1}$ has trivial Jacobian,
$J_{\Proj^{1}, \mathfrak{m}}$ is an affine algebraic group times $\Z$.
The groups $J_{\Proj^{1}, \mathfrak{m}}$, the schemes $\Proj^{1} - \Supp{\mathfrak{m}}$
and the morphisms $\Phi_{\mathfrak{m}} \colon \Proj^{1} - \Supp{\mathfrak{m}} \to J_{\Proj^{1}, \mathfrak{m}}$
form projective systems with respect to $\mathfrak{m}$.
Let $J_{\Proj^{1}}$ (resp.\ $\Phi \colon \Spec k(T) \to J_{\Proj^{1}}$) be
the projective limit of the $J_{\Proj^{1}, \mathfrak{m}}$ (resp.\ the $\Phi_{\mathfrak{m}}$),
the limit being taken over all effective divisors $\mathfrak{m}$ on $\Proj^{1}$.
They are naturally regarded as sheaves on $\pfpqc{k}$
and morphisms of sheaves on $\pfpqc{k}$.
The inclusion $k(T) \subset k((T)) = K$ identifies $K$
as the fraction field of the completion of the local ring of $\Proj^{1}$ at $T = 0$.
This identification gives a natural morphism $\Spec K_{p} \to \Spec k(T)_{p}$
and a natural homomorphism $\alg{K}^{\times} \to J_{\Proj^{1}}$.
Consider the following (non-commutative) diagram:
	\[
		\begin{CD}
				\Spec K_{p}
			@>> \varphi >
				\alg{K}^{\times}
			\\
			@VVV
			@VVV
			\\
				\Spec k(T)_{p}
			@> \Phi >>
				J_{\Proj^{1}}.
		\end{CD}
	\]
This diagram gives two morphisms $\Spec K_{p} \to J_{\Proj^{1}}$,
or equivalently, two $K_{p}$-rational points of $J_{\Proj^{1}}$.
Under the natural identification
	$
			J_{\Proj^{1}}
		=
			(
				\prod_{a \in \Proj^{1}(k) - \{0\}}
					\alg{U}_{a} \times \alg{K}^{\times}
			) / \Gm
	$
with $\alg{U}_{a}$ defined by
	$
			\alg{U}_{a}(R)
		=
			R[[\alg{T}^{-1} - a^{-1}]]^{\times}
	$,
the difference of these two rational points is given by
	$
			((1 - T \alg{T}^{-1})_{a \in \Proj^{1}(k) - \{0\}}, 1)
	$,
which is in $J_{\Proj^{1}}((\maxid_{K})_{p})$.
Hence this diagram induces, on fundamental groups, a commutative diagram
	\[
		\begin{CD}
				\Gal(K^{\ab} / K)
			@>> \eta >
				\pi_{1}^{k}(\alg{K}^{\times})
			\\
			@VVV
			@VVV
			\\
				\Gal(k(T)^{\ab} / k(T))
			@>>>
				\pi_{1}^{k}(J(\Proj^{1}))
		\end{CD}
	\]
by the same argument as in the proof of Proposition \ref{Prop:Indep}.
The top horizontal map $\eta$ is an isomorphism
by Assertion \ref{Ass:Main:Isom} of Theorem \ref{Thm:Main}.
The bottom horizontal map is an isomorphism as well by \cite{Ser88}.
Regarding this diagram,
Serre proved the global version of the equality $\eta^{-1} = \theta$
and the local-global compatibility
(\cite[\S 5.2, Prop.\ 2 and \S 5.3, Lem.\ 2]{Ser61};
actually Serre's Proposition 2, which says $\psi^{-1} = \theta$ in his notation,
should be corrected as $- \psi^{-1} = \theta$).
This implies the local equality $\eta^{-1} = \theta$.


\section{A two-dimensional analog}
\label{Sect:TwoDim}

In this section, we discuss an analog of the above theory
for a field of the form $K = k((S))((T))$.
We denote by $K_{2}$ the functor of the second algebraic $K$-group (\cite{Bas73}).
For each perfect $k$-algebra $R$,
we have the abelian group $K_{2}(R[[\alg{S}, \alg{T}]])$.
This gives a group functor that we denote by $K_{2}[[\alg{S}, \alg{T}]]$.
The $K_{p}$-rational point
	\[
			\{- S + \alg{S}, - T + \alg{T}\}
		\in
			K_{2}[[\alg{S}, \alg{T}]](K_{p})
		=
			K_{2} \bigl(
				k((S))((T))_{p}
				[[\alg{S}, \alg{T}]]
			\bigr)
	\]
gives a morphism
	$
			\varphi^{0}
		\colon
			\Spec K_{p}
		\to
			K_{2}[[\alg{S}, \alg{T}]]
	$,
where $\{\cdot, \cdot\}$ denotes the symbol map.
This is an analog of the map $\Spec k((T))_{p} \to \alg{U}_{k((T))}$
we have defined and studied before this section.

Now assume $k$ is the finite field $\F_{q}$ with $q$ elements.
Define a functor $X$ by the following cartesian diagram:
	\[
		\begin{CD}
				X
			@>>>
				K_{2}[[\alg{S}, \alg{T}]]
			\\
			@VVV
			@VV F - 1 V
			\\
				\Spec K_{p}
			@> \varphi^{0} >>
				K_{2}[[\alg{S}, \alg{T}]],
		\end{CD}
	\]
where $F$ is the $q$-th power Frobenius morphism over $k$.
Then we expect that $K_{2}[[\alg{S}, \alg{T}]]$ can be viewed
as a sort of an ``algebraic group over $k$''
and the equation $(F - 1) x = \{- S + \alg{S}, - T + \alg{T}\}$
gives a two-dimensional analog of Lubin-Tate theory
so that $X$ is the $\Spec$ of the perfect closure of
a large totally ramified abelian extension of $K$
(compare with Proposition \ref{Prop:LT}).

Although this hope is hard to be made precise,
we can formulate and prove a rigorous partial result
(Proposition \ref{Prop:TwoDim} below).
For each perfect $k$-algebra $R$,
we have the space of $1$-forms
$\Omega_{R[[\alg{S}, \alg{T}]] / R}^{1}$
and the space of $2$-forms
$\Omega_{R[[\alg{S}, \alg{T}]] / R}^{2}$.
We denote these functors by
$\Omega^{1}_{[[\alg{S}, \alg{T}]]}$, $\Omega^{2}_{[[\alg{S}, \alg{T}]]}$,
respectively.
We have
	\begin{gather*}
				\Omega^{1}_{[[\alg{S}, \alg{T}]]}
			=
				\prod_{i, j \ge 0}
						\Ga \alg{S}^{i} \alg{T}^{j} d \alg{S}
					\times
						\Ga \alg{S}^{i} \alg{T}^{j} d \alg{T},
		\\
				\Omega^{2}_{[[\alg{S}, \alg{T}]]}
			=
				\prod_{i, j \ge 1}
						\Ga \alg{S}^{i} \alg{T}^{j}
						\dlog \alg{S} \wedge \dlog \alg{T},
	\end{gather*}
where $\dlog \alg{S} = d \alg{S} / \alg{S}$,
$\dlog \alg{T} = d \alg{T} / \alg{T}$.
Consider the following sequence:
	\[
			0
		\longrightarrow
			K_{2}[[\alg{S}, \alg{T}]] / p K_{2}[[\alg{S}, \alg{T}]]
		\overset{\dlog}{\longrightarrow}
			\Omega^{2}_{[[\alg{S}, \alg{T}]]}
		\overset{C^{-1} - 1}{\longrightarrow}
				\Omega^{2}_{[[\alg{S}, \alg{T}]]}
			/
				d \Omega^{1}_{[[\alg{S}, \alg{T}]]}.
	\]
Here $\dlog$ is the dlog map of algebraic $K$-theory and
$C^{-1}$ is the inverse Cartier operator.
We do not claim the exactness of this sequence;
however, it should be compared with Bloch-Kato-Gabber's theorem
(\cite[Chapter 2, \S 2.4, Theorem 5]{FK00}).
We study $\Ker(C^{-1} - 1)$ instead of
$K_{2}[[\alg{S}, \alg{T}]] / p K_{2}[[\alg{S}, \alg{T}]]$.
We have a commutative diagram
	\[
		\begin{CD}
				X
			@>>>
				K_{2}[[\alg{S}, \alg{T}]]
			@>> \dlog >
				\Ker(C^{-1} - 1)
			\\
			@VVV
			@VV F - 1 V
			@VV F - 1 V
			\\
				\Spec K_{p}
			@> \varphi^{0} >>
				K_{2}[[\alg{S}, \alg{T}]]
			@> \dlog >>
				\Ker(C^{-1} - 1).
		\end{CD}
	\]
We put $\varphi^{1} = {\dlog} \circ \varphi^{0}$.

\begin{prop}
	\label{Prop:TwoDim} \mbox{}
	\begin{enumerate}
		\item \label{Ass:TwoDim:Cartesian}
			There is a cartesian diagram
				\[
					\begin{CD}
							\Spec L_{p}
						@>>>
							\Ker(C^{-1} - 1)
						\\
						@VVV
						@VV F - 1 V
						\\
							\Spec K_{p}
						@> \varphi^{1} >>
							\Ker(C^{-1} - 1).
					\end{CD}
				\]
			Here we put
				$
						L
					:=
						K[x_{i j} \,|\, i, j \ge 1,\, p \nmid \gcd(i, j)]
					/
						(x_{i j}^{q} - x_{i j} - S^{- i} T^{- j})
				$.
		\item \label{Ass:TwoDim:Field}
			The ring $L$ is actually a field.
			The extension $L / K$ is an abelian (infinite) extension
			with Galois group isomorphic to the group
			of $k$-rational points of $\Ker(C^{-1} - 1)$
			(see \eqref{Eq:KerCartier} below for the structure of this group).
	\end{enumerate}
\end{prop}

\begin{proof}
	\ref{Ass:TwoDim:Cartesian}.
	First we show that the projection gives an isomorphism
		\begin{equation} \label{Eq:KerCartier}
				\Ker(C^{-1} - 1)
			\isomto
				\prod_{i, j \ge 1,\, p \nmid \gcd(i, j)}
					\Ga \alg{S}^{i} \alg{T}^{j}
					\dlog \alg{S} \wedge \dlog \alg{T}.
		\end{equation}
	We have
		$
					\Omega^{2}_{[[\alg{S}, \alg{T}]]}
				/
					d \Omega^{1}_{[[\alg{S}, \alg{T}]]}
			=
				\prod_{i, j \ge 1}
					\Ga \alg{S}^{p i} \alg{T}^{p j}
					\dlog \alg{S} \wedge \dlog \alg{T}
		$.
	The map $C^{-1} - 1$ sends an element
		$
				\sum_{i, j \ge 1}
					a_{i j} \alg{S}^{i} \alg{T}^{j}
					\dlog \alg{S} \wedge \dlog \alg{T}
		$
	of $\Omega^{2}_{[[\alg{S}, \alg{T}]]}$
	to
		\begin{align*}
			&
					\sum_{i, j \ge 1}
						(a_{i j}^{p} \alg{S}^{p i} \alg{T}^{p j}
						- a_{i j} \alg{S}^{i} \alg{T}^{j})
						\dlog \alg{S} \wedge \dlog \alg{T}
			\\
			&	=
					\sum_{i, j \ge 1}
						(a_{i j}^{p} - a_{p i \, p j}) \alg{S}^{p i} \alg{T}^{p j}
						\dlog \alg{S} \wedge \dlog \alg{T}
			\\
			&	\in
						\Omega^{2}_{[[\alg{S}, \alg{T}]]}
					/
						d \Omega^{1}_{[[\alg{S}, \alg{T}]]}
				=
					\prod_{i, j \ge 1}
						\Ga \alg{S}^{p i} \alg{T}^{p j}
						\dlog \alg{S} \wedge \dlog \alg{T}.
		\end{align*}
	Hence 
		$
				\sum_{i, j \ge 1}
					a_{i j} \alg{S}^{i} \alg{T}^{j}
					\dlog \alg{S} \wedge \dlog \alg{T}
		$
	is in $\Ker(C^{-1} - 1)$ if and only if
	$a_{i j}^{p} = a_{p i \, p j}$ for any $i, j \ge 1$.
	This gives the required description of $\Ker(C^{-1} - 1)$.
	
	The map $\varphi^{1} \colon \Spec K_{p} \to \Ker(C^{-1} - 1)$
	corresponds to the rational point
		\begin{align*}
					\dlog \{-S + \alg{S}, -T + \alg{T}\}
			&	=
						\frac{d(-S + \alg{S})}{-S + \alg{S}}
					\wedge
						\frac{d(-T + \alg{T})}{-T + \alg{T}}
			\\
			&	=
					\sum_{i, j \ge 1}
						S^{-i} T^{-j}
						\alg{S}^{i} \alg{T}^{j}
						\dlog \alg{S} \wedge \dlog \alg{T},
		\end{align*}
	which corresponds to
		$
			\sum_{i, j \ge 1,\, p \nmid \gcd(i, j)}
				S^{-i} T^{-j}
				\alg{S}^{i} \alg{T}^{j}
				\dlog \alg{S} \wedge \dlog \alg{T}
		$
	via the above isomorphism.
	If
		$
			\sum_{i, j \ge 1,\, p \nmid \gcd(i, j)}
				x_{i j}
				\alg{S}^{i} \alg{T}^{j}
				\dlog \alg{S} \wedge \dlog \alg{T}
		$
	lies in the geometric fiber of $F - 1$ over this rational point,
	it must satisfy
		\begin{align*}
			&
					(F - 1)
					\sum_{i, j \ge 1,\, p \nmid \gcd(i, j)}
						x_{i j}
						\alg{S}^{i} \alg{T}^{j}
						\dlog \alg{S} \wedge \dlog \alg{T}
			\\
			&	=
					\sum_{i, j \ge 1,\, p \nmid \gcd(i, j)}
						(x_{i j}^{q} - x_{i j})
						\alg{S}^{i} \alg{T}^{j}
						\dlog \alg{S} \wedge \dlog \alg{T}
			\\
			&	=
					\sum_{i, j \ge 1,\, p \nmid \gcd(i, j)}
						S^{-i} T^{-j}
						\alg{S}^{i} \alg{T}^{j}
						\dlog \alg{S} \wedge \dlog \alg{T}
		\end{align*}
	This proves Assertion \ref{Ass:TwoDim:Cartesian}.
	
	\ref{Ass:TwoDim:Field}.
	By Artin-Schreier theory, we have an isomorphism
		$
				K / \wp_{q}(K)
			\isomto
				\Hom(\Gal(K^{\sep} / K), \F_{q})
		$
	which sends an element $f \in K / \wp_{q} K$ to the character
	$\sigma \to (\sigma - 1) \wp_{q}^{-1}(f)$,
	where $\wp_{q}(x) = x^{q} - x$.
	Consider the (non-commutative) $\F_{q}$-algebra
	$\sgpalg{\F_{q}}{\Gal(\F_{q} / \F_{p})}$
	with $\F_{q}$-basis $\Gal(\F_{q} / \F_{p})$
	and relations $\sigma a = \sigma(a) \sigma$
	for each $a \in \F_{q}$ and $\sigma \in \Gal(\F_{q} / \F_{p})$.
	It naturally acts on $\F_{q}$ and hence on $\Hom(\Gal(K^{\sep} / K), \F_{q})$.
	The action of the $p$-th power Frobenius $\in \Gal(\F_{q} / \F_{p})$
	corresponds to the $p$-th power map on $K / \wp_{q}(K)$
	via the isomorphism
		$
				K / \wp_{q}(K)
			\isomto
				\Hom(\Gal(K^{\sep} / K), \F_{q})
		$.
	By the following lemma,
	the proof for $L$ being a field is reduced to
	showing that the elements $S^{-i} T^{-j}$ for $i, j \ge 1$, $p \nmid \gcd(i, j)$
	generate a free $\sgpalg{\F_{q}}{\Gal(\F_{q} / \F_{p})}$-submodule
	of $K / \wp_{q}(K)$.
	
	\begin{lemm}
		Let $\chi_{1}, \dots, \chi_{n}$ be $n$ elements of $\Hom(\Gal(K^{\sep} / K), \F_{q})$
		and let $K_{i}$ correspond $\Ker(\chi_{i})$ via Galois theory, $1 \le i \le n$.
		Assume $\chi_{1}, \dots, \chi_{n}$ generate
		a free $\sgpalg{\F_{q}}{\Gal(\F_{q} / \F_{p})}$-submodule
		of $\Hom(\Gal(K^{\sep} / K), \F_{q})$.
		Then the fields $K_{1}, \dots, K_{n}$ have degree $q$
		and are linearly disjoint.
	\end{lemm}
	
	\begin{proof}[Proof of the lemma]
		The natural map
			$
					\sgpalg{\F_{q}}{\Gal(\F_{q} / \F_{p})}
				\to
					\End_{\text{$\F_{p}$-linear}}(\F_{q})
			$
		is an $\F_{q}$-algebra isomorphism
		since it is injective by the linear independence of field automorphisms
		and they have the same dimension over $\F_{q}$.
		Let $\psi = (\psi_{1}, \dots, \psi_{n})$ be an $\F_{p}$-linear map
		$\F_{q}^{\dsum n} \to \F_{q}$ which vanishes on the image of
		$\Dsum_{i} \chi_{i} \colon \Gal(K^{\sep} / K) \to \F_{q}^{\dsum n}$.
		Then $\sum_{i} \psi_{i} \circ \chi_{i} = 0$.
		By this and the assumption on $\chi_{i}$,
		we have $\psi_{i} = 0$ for any $i$.
		Hence $\psi = 0$.
		This implies the surjectivity of
		$\Dsum_{i} \chi_{i} \colon \Gal(K^{\sep} / K) \to \F_{q}^{\dsum n}$.
		This proves the lemma.
	\end{proof}
	
	The additive group $K$ can be written as the direct sum
	of $V := \Dsum_{i, j \ge 1} \F_{q} S^{-i} T^{-j}$
	and $\{\sum_{\max\{i, j\} \ge 0} a_{i j} S^{i} T^{j} \in K\}$,
	both being stable under the action of $\wp_{q}$.
	The space $V / \wp_{q}(V)$ injects into $K / \wp_{q}(K)$
	and has $S^{-i} T^{-j}$ for $i, j \ge 1$, $q \nmid \gcd(i, j)$ as an $\F_{q}$-basis.
	Hence $V / \wp_{q}(V)$ is a free $\sgpalg{\F_{q}}{\Gal(\F_{q} / \F_{p})}$-module
	with basis $S^{-i} T^{-j}$ for $i, j \ge 1$, $p \nmid \gcd(i, j)$.
	This implies that $L$ is a field.
	
	The latter half of Assertion \ref{Ass:TwoDim:Field} follows from
	the cartesian diagram of Assertion \ref{Ass:TwoDim:Cartesian}.
\end{proof}

\begin{rema}
	By a similar argument to the proof of Assertion \ref{Ass:TwoDim:Field},
	we have
		\begin{align*}
			&		\Hom(\Gal(K^{\sep} / K), \F_{q})
				\cong
					K / \wp_{q} K
			\\
			&	\cong
						\F_{q}
					\dsum
						\Dsum_{q \nmid i \ge 1} \F_{q} S^{- i}
					\dsum
						\Dsum_{q \nmid j \ge 1} \F_{q} T^{- j}
					\dsum
			\\
			&	\quad
						\Dsum_{i \ge 1} \prod_{q \nmid j \ge 1}
							\F_{q} S^{i} T^{- j}
					\dsum
						\Dsum_{q \nmid i \ge 1} \prod_{j \ge 1}
							\F_{q} S^{i} T^{- j}
					\dsum
						\Dsum_{i, j \ge 1,\, q \nmid \gcd(i, j)}
							\F_{q} S^{- i} T^{- j}.
		\end{align*}
	The last summand $\Dsum \F_{q} S^{- i} T^{- j}$ corresponds to
	the abelian extension $L$ in the proposition.
\end{rema}


\section{$\D$-modules on two-dimensional $K$}
\label{Sect:D-mod}

In this section, we give an analog of the result of \cite[\S 2.6]{BBDE} on $\D$-modules,
for a field of the form $K = k((S))((T))$ with $k$ characteristic zero.
This is also regarded as an analog of what we have done in the previous section.

We formulate our result (Proposition \ref{Prop:D-mod} below).
For a $k$-algebra $R$,
we define $\Omega_{[[\alg{S}, \alg{T}]]}^{2}(R)$ to be the projective limit of
the space of $2$-forms of $R[\alg{S}, \alg{T}] / (\alg{S}, \alg{T})^{n}$ over $R$,
the limit being taken for $n \ge 0$.
Let $\varphi^{1} \colon \Spec K \to \Omega_{[[\alg{S}, \alg{T}]]}^{2}$
be the $k$-morphism corresponding to the rational point
$\dlog \{- S + \alg{S}, - T + \alg{T}\} = \dlog(- S + \alg{S}) \wedge \dlog(- T + \alg{T})$.
We say that a $\D$-module $M$ of $\Order$-rank one on $\Omega_{[[\alg{S}, \alg{T}]]}^{2}$ is
compatible with the group structure if
$\mu^{\ast} M \cong \pr_{1}^{\ast} M \tensor \pr_{2}^{\ast} M$,
where
	$
			\mu
		\colon
			\Omega_{[[\alg{S}, \alg{T}]]}^{2}
		\times
			\Omega_{[[\alg{S}, \alg{T}]]}^{2}
		\to
			\Omega_{[[\alg{S}, \alg{T}]]}^{2}
	$
is the addition map and $\pr_{i}$ is the $i$-th projection for each $i = 1, 2$.
Let $\mathcal{C}$ be the category of $\D$-modules of $\Order$-rank one
on $\Omega_{[[\alg{S}, \alg{T}]]}^{2}$ that are compatible with the group structure
and let $\mathcal{C}'$ be the category of $\D$-modules of $\Order$-rank one on $\Spec K$.
The categories $\mathcal{C}$ and $\mathcal{C}'$ are
$k$-linear Picard categories under the tensor operation.
The endomorphism ring of any object of $\mathcal{C}$ or $\mathcal{C}'$ is equal to $k$.
Let $\Isom(\mathcal{C})$ (resp.\ $\Isom(\mathcal{C}')$) be
the abelian group of isomorphism classes of objects of $\mathcal{C}$ (resp.\ $\mathcal{C}'$).
Since $\Omega_{[[\alg{S}, \alg{T}]]}^{2}$ and $\Spec K$ are the $\Spec$'s of unique factorization domains,
line bundles on them can be trivialized.
We can associate connection forms
to objects of $C$ and $C'$ whose underlying line bundles are trivial.
Thus $\Isom(\mathcal{C})$ (resp.\ $\Isom(\mathcal{C}')$) can identified with
a subquotient (as an abelian group) of the space of $1$-forms
on the $k$-scheme $\Omega_{[[\alg{S}, \alg{T}]]}^{2}$ (resp.\ $\Spec K$).

\begin{prop}
	\label{Prop:D-mod} \mbox{}
	\begin{enumerate}
		\item \label{Ass:D-mod:Emb}
			The pullback by $\varphi^{1} \colon \Spec K \to \Omega_{[[\alg{S}, \alg{T}]]}^{2}$
			gives a fully faithful embedding $\mathcal{C} \into \mathcal{C}'$
			of $k$-linear Picard categories.
		\item \label{Ass:D-mod:Isom(C')}
			Under the above identification of $\mathcal{C}'$, we have
				\begin{align*}
							\Isom(\mathcal{C}')
					&	\cong
							(k((T)) / \Z) \dlog S \dsum (k((S)) / \Z) \dlog T
					\\
					&	\qquad \dsum
							d \Bigl (
									S^{-1} k[S^{-1}] \dsum T^{-1} k[T^{-1}]
					\\
					&	\qquad \qquad
								\dsum
									T^{-1} S k[[S]][T^{-1}]
								\dsum
									T^{-1} S^{-1} k[S^{-1}, T^{-1}]
							\Bigr ).
				\end{align*}
		\item \label{Ass:D-mod:Isom(C)}
			The image of $\Isom(\mathcal{C})$ in $\Isom(\mathcal{C}')$
			by the pullback functor by $\varphi^{1}$
			corresponds to the last summand
			$d (T^{-1} S^{-1} k[S^{-1}, T^{-1}])$
			via the isomorphism of Assertion \ref{Ass:D-mod:Isom(C')}.
	\end{enumerate}
\end{prop}

\begin{proof}
	\ref{Ass:D-mod:Isom(C')}.
	The group $\Isom(\mathcal{C}')$ is identified with
	$\Omega_{K / k}^{1, d = 0} / \dlog K^{\times}$.
	A straightforward calculation shows that
	$\Omega_{K / k}^{1, d = 0} / \dlog K^{\times}$ is equal to
	the right-hand side of Assertion \ref{Ass:D-mod:Isom(C')}.
	
	\ref{Ass:D-mod:Isom(C)}.
	The group $\Omega_{[[\alg{S}, \alg{T}]]}^{2}$
	is the product of copies of $\Ga$ labeled by $n, m \ge 1$
	with coordinates $z_{n m} := \alg{S}^{n} \alg{T}^{m} \dlog \alg{S} \wedge \dlog \alg{T}$.
	Hence a connection form of any object $M$ of $\mathcal{C}$
	is of the form $\sum_{n, m \ge 1} a_{n m} d z_{n m}$
	for some $a_{n m} \in k$ that are uniquely determined by the isomorphism class of $M$.
	Therefore $\Isom(\mathcal{C}) \cong \Dsum_{n, m \ge 1} k d_{n m}$.
	Since
		$
				\dlog \{- S + \alg{S}, - T + \alg{T}\}
			=
				\sum_{n, m \ge 1}
					S^{- n} T^{- m} z_{n m}
		$,
	the pullback of a connection form $\sum_{n, m \ge 1} a_{n m} d z_{n m}$
	with $a_{n m} \in k$ is
		$
				\sum_{n, m \ge 1}
					a_{n m} d (S^{- n} T^{- m})
			=
				d (
					\sum_{n, m \ge 1}
					a_{n m} S^{- n} T^{- m}
				)
		$.
	Thus we get the result.
	
	\ref{Ass:D-mod:Emb}.
	The above proof of Assertion \ref{Ass:D-mod:Isom(C)} also shows that
	the homomorphism $\Isom(\mathcal{C}) \to \Isom(\mathcal{C}')$
	induced by the functor $\mathcal{C} \to \mathcal{C}'$ is injective.
	The result follows from this.
	
\end{proof}


\bibliographystyle{smfalpha}
\bibliography{Suzuki_LCFT.bib}

\end{document}